%% file: main.tex
\documentclass{article}
\usepackage{icomp2026_conference,times}
\usepackage{amsmath,amssymb,amsthm,mathtools,booktabs,tabularx,array,microtype,enumitem}
\usepackage[most]{tcolorbox}
\usepackage{tikz}
\usetikzlibrary{arrows.meta,positioning,shapes.geometric}
\usepackage{hyperref}
\usepackage{url}
\def\doi#1{\href{https://doi.org/#1}{doi:\nolinkurl{#1}}}
\hypersetup{hidelinks,pdfauthor={Amir Valiakhmetov, Ewsey Obzherin, Alexey Chernov, Roland Hildebrand},pdftitle={Four-Dimensional Semi-Homogeneous Convex Cones: A Complete Classification}}
\allowdisplaybreaks

\title{Four-Dimensional Semi-Homogeneous Convex Cones: A Complete Classification}
\author{
Amir Valiakhmetov
\And
Ewsey Obzherin\thanks{Corresponding author: obzherin.ewsey@ya.ru}
\And
Alexey Chernov
\And
Roland Hildebrand
}

\newtheorem{theorem}{Theorem}
\newtheorem{proposition}[theorem]{Proposition}
\newtheorem{lemma}[theorem]{Lemma}

\theoremstyle{definition}
\newtheorem{definition}{Definition}

\tcbset{keybox/.style={
  enhanced,breakable,colback=white,colframe=black,boxrule=0.45pt,arc=1pt,
  left=4pt,right=4pt,top=3pt,bottom=3pt,before skip=6pt,after skip=6pt}}
\tcolorboxenvironment{theorem}{keybox}
\tcolorboxenvironment{definition}{keybox}

\newcommand{\R}{\mathbb R}
\newcommand{\Rp}{\mathbb R_{+}}
\newcommand{\Aut}{\operatorname{Aut}}
\newcommand{\aut}{\mathfrak{aut}}
\newcommand{\cl}{\operatorname{cl}}
\newcommand{\diag}{\operatorname{diag}}
\newcommand{\Span}{\operatorname{span}}
\newcommand{\rank}{\operatorname{rank}}
\newcommand{\tr}{\operatorname{tr}}
\newcommand{\e}{\mathrm e}

\newcommand{\cL}{\mathcal C_L}
\newcommand{\EndPairs}{\mathcal E}

\setlist[itemize]{leftmargin=1.45em,itemsep=1pt,topsep=2pt}
\setlist[enumerate]{leftmargin=1.65em,itemsep=1.5pt,topsep=2pt}

\providecommand{\keywords}[1]{%
  \par\noindent\textbf{Keywords: }\begingroup\def\and{\unskip;\space}#1\endgroup\par
}

\begin{document}
\maketitle

\begin{abstract}
We classify all four-dimensional semi-homogeneous convex cones, i.e., closed, pointed, full-dimensional cones whose connected automorphism group has principal orbits of codimension at most one. On such cones, the Monge-Amp\`ere equation defining the canonical barrier reduces to an ordinary differential equation. This makes this barrier computationally accessible. The semi-homogeneous cones can then be used in modeling by conic programs. An obvious subset of such cones are the homogeneous cones. The key point allowing the classification is that for all other cones, the automorphism group has a three-dimensional subgroup acting locally faithfully on the principal orbits. This subgroup must have a Lie algebra which is either abelian or a direct product of a 1-dimensional algebra and a two-dimensional algebra generated by $A,B$ satisfying $[A,B]=B$. In both cases we classify all 4-dimensional representations up to isomorphism, and for each representation all invariant convex cones. For every cone an orbit invariant is computed.
\end{abstract}


\input{sections/section_01_introduction}
\input{sections/section_02_semi_homogeneous_cones_and_notation}
\input{sections/section_03_the_classification}
\input{sections/section_04_reduction_to_two_matrix_equations}
\input{sections/section_05_the_abelian_case}
\input{sections/section_06_the_nonabelian_case}
\input{sections/section_07_conclusion}

\clearpage
\bibliographystyle{icomp2026_conference}
\bibliography{references}

\clearpage
\appendix
\input{sections/appendix_structure}
\input{sections/appendix_abelian}
\input{sections/appendix_affine}

\input{sections/appendix_verification}

\end{document}

%% file: sections/section_01_introduction.tex
\section{Introduction}

For every closed pointed convex cone $K\subset\R^n$ with nonempty interior $K^{\circ}$, the canonical potential, or canonical barrier, is the unique solution $F: K^{\circ} \to \mathbb R$ of the Monge-Amp\`ere equation
\begin{equation}\label{eq:MA-main}
 \det D^2F=\e^{2F},\qquad D^2F\succ0,\qquad F|_{\partial K}=+\infty.
\end{equation}
This solution automatically satisfies the logarithmic homogeneity condition
\[ F(tx)=F(x)-n\log t.
\]
More generally, if $g \in GL(n,\mathbb R)$ is a linear automorphism of the cone $K$, then \citep{HildebrandCanonical2014}
\begin{equation} \label{automorphismSymmetry} 
F(gx) = F(x) - \log|\det\,g|.
\end{equation}

The level sets of $F$ are the complete hyperbolic affine spheres asymptotic to $\partial K$~\citep{Calabi1972,ChengYau1980,ChengYau1982,ChengYau1986,LiSimonZhao1993,NomizuSasaki1994,Loftin2001,Loftin2002,Sasaki1985,HildebrandCanonical2014,Fox2015}. Existence of the canonical barrier is assured for general cones, but explicit formulas are rare. If the automorphism group $\Aut\,K$ has orbits of codimension one, i.e., the cone is \textit{semi-homogeneous}, then relation \eqref{automorphismSymmetry} allows one to reduce the partial differential equation (PDE) \eqref{eq:MA-main} on the domain $K^{\circ} \subset \mathbb R^n$ to an ordinary differential equation (ODE) on the one-dimensional orbit space. In dimension three all semi-homogeneous cones were classified and the resulting ODEs solved in~\citep{HildebrandAffine2014}. The ODEs turned generically out to be solvable by elliptic integrals. The semi-homogeneous affine-sphere families foliating the interior of the cones were explicitly parameterized in~\citep{LinWang2016}.

This reduction of PDE \eqref{eq:MA-main} to an ODE can be put to use in conic optimization. Computable logarithmically homogeneous self-concordant barriers are necessary for polynomial-time interior-point methods to work~\citep{NesterovNemirovskii1994}. Such barriers are readily available on symmetric cones. Universal and canonical barriers exist on arbitrary cones~\citep{Guler1996,HildebrandCanonical2014,HildebrandProjective2022,BubeckEldan2019,LeeYue2021} but are in general not efficiently computable. For power, exponential, entropy, and other nonsymmetric cones specialized barrier oracles have been developed~\citep{Nesterov2012,SkajaaYe2015,PappYildiz2022,CoeyKapelevichVielma2022,CoeyKapelevichVielma2023,PappVarga2025}. Generalized power cones provide a recent example where automorphisms and facial structure can be described explicitly~\citep{LinEtAl2024}. 

Here we provide a complete classification of four-dimensional semi-homogeneous cones. On these cones the canonical barrier is more easily computable due to the reduction of PDE \eqref{eq:MA-main} to an ODE. Among the semi-homogeneous cones we find the already known generalized power cones, but also new families of cones which can now be considered as building blocks for modeling optimization problems by conic programs.

The classification of cones is carried out via a classification of the Lie algebra $\aut(K)$ of the automorphism group or a subalgebra thereof. Homotheties are always part of $\Aut\,K$ and correspond to multiples of the identity matrix in $\aut(K)$. After they are factored out, the subgroup preserving the level surfaces of the canonical barrier remains. The corresponding subalgebra of $\aut(K)$ contains exactly the traceless matrices.

We shall show that the following cases are possible:
\begin{itemize}
    \item $K$ is homogeneous
    \item the principal orbits (those with largest dimension, which is also the generic dimension) of $\Aut(K)$ have dimension 3 and there exists a connected 3-dimensional subgroup $G$ of $\Aut(K)$, which includes the homotheties, with dimension of its principal orbits also equal 3
\end{itemize}

The homogeneous cones are $\mathbb R_+^4$, $L_4$, and $L_3 \times \mathbb R_+$, up to linear isomorphism \citep{KaneyukiTsuji74,Vinberg1963}. Here $L_n$ is the $n$-dimensional Lorentz or second-order cone. In the non-homogeneous case the interior of the cone splits into orbits of generic codimension 1 (Lemma \ref{lem:suffCondition}). The family of orbits can be parameterized by a continuous scalar function, the orbit invariant. The interior of the cone must then be the union of orbits with invariants taking values in some interval. It will be shown below in Theorem \ref{thm:reduction-main} and Lemma \ref{lem:suffCondition} that the two cases above are exhaustive.

In the non-homogeneous case, the Lie subalgebra $\mathfrak g_0 \subset \mathfrak g$ of trace-less matrices in the Lie algebra of $G$ is 2-dimensional. It is then either abelian or there exist basis elements $A,B$ of $\mathfrak g_0$ satisfying  $[A,B]=B$. In the commuting case the simultaneous generalized-weight decomposition gives five dimension patterns. For each pattern, we derive a canonical basis of $\mathfrak g_0$, compute the complete scalar orbit invariant, establish cone convexity either through explicit inequalities or via the Hessian of a boundary function, and identify all invariant proper cones. In the second case, the relation $[A,B]=B$ is explicitly solved and the invariant cones are computed for each emerging subcase.

The paper is organized as follows. Section~\ref{sec:2} fixes the terminology and the notation. In Section~\ref{sec:3} we state the classification and give tables with a summary of the results. Section~\ref{sec:4} establishes that the homogeneous and non-homogeneous cases described above are exhaustive. Sections~\ref{sec:5} and~\ref{sec:6} contain the classification for the abelian and nonabelian subcases of the non-homogeneous case, respectively. The appendices contain the detailed proofs, and the treatment of those cases not leading to invariant convex cones.

%% file: sections/section_02_semi_homogeneous_cones_and_notation.tex
\section{Semi-homogeneous cones and notation} \label{sec:2}

Throughout the paper we consider cones in $\mathbb R^4$.

\begin{definition}[Automorphisms]\label{def:actions}
Let $K \subset \mathbb R^4$ be a conic set. We define its automorphism group and the Lie algebra thereof by
\[
\begin{aligned}
 \Aut(K)&=\{g\in\mathrm{GL}_4(\R):gK=K\},\\
 \aut(K)&=\{A\in\mathfrak{gl}_4(\R):\e^{tA}K=K
          \text{ for all }t\in\R\}.
\end{aligned}
\]
Positive homotheties belong to $\Aut(K)$, hence $I\in\aut(K)$.
\end{definition}

\begin{definition}[Proper and (semi-)homogeneous cones]\label{def:semihomogeneous}
A cone $K\subset\R^4$ is \emph{proper} if it is closed, convex, pointed, and has nonempty interior. It is \emph{semi-homogeneous} if a principal orbit of $\Aut(K)$ has codimension at most one in $K^\circ$. It is homogeneous if its interior is an orbit of its automorphism group.
\end{definition}

If the principal orbits have full dimension 4, then $K$ is homogeneous by virtue of Lemma \ref{lem:suffCondition} below. The non-homogeneous semi-homogeneous cones have three-dimensional principal orbits. As mentioned in the introduction, for every $g\in\Aut(K)$ relation \eqref{automorphismSymmetry} is valid \citep{HildebrandCanonical2014}.

Introduce coordinates $(x,y,z,w)$ in $\mathbb R^4$. Put
\begin{equation}\label{eq:lorentz-main}
 \Delta=2xz-y^2,
 \qquad
 \cL=\{(x,y,z):x\ge0,\ z\ge0,\ \Delta\ge0\},
 \qquad
 H_\gamma=\Delta^{\gamma/2}z^{1-\gamma}.
\end{equation}
Note that $H_{\gamma}$ is a 1-homogeneous function, and $\cL$ is isomorphic to the Lorentz cone $L_3$.

We shall work with open intervals $(-a,b) \not= \mathbb R$ of the real line containing zero in their closure. We consider two intervals as equivalent if one can be obtained from the other by multiplication by a non-zero real number. Then canonical representatives of the resulting equivalence classes are given by pairs of numbers $(a,b)$ in the set
\begin{equation}\label{eq:endpairs-main}
 \EndPairs=\{(a,1):0\leq a\le1\}\cup\{(0,+\infty),(1,+\infty)\}
\end{equation}
defining intervals $(-a,1),(0,+\infty),(-1,\infty)$, respectively. %Reflection about 0 and positive rescaling of the real line $\mathbb R$ will reduce every nontrivial interval with endpoints of opposite signs, or with one finite endpoint and not strictly separated from 0, to exactly one pair in~\eqref{eq:endpairs-main}.

Some cone families are described by constraints containing expressions like $x^{1/p}$, $p \in (1,+\infty)$. Let $q=p/(p-1)$ and define the three-dimensional cones
\begin{align}
 K_{\exp}^{3}&=\cl\{(x,y,z):z>0,\ y>z\e^{x/z}\},\label{eq:Kexp3-main}\\
 K_{\alpha,\beta}^{p,3}&=\cl\{(x,y,z):x,z>0,
 -\alpha x^{1/p}z^{1/q}<y<\beta x^{1/p}z^{1/q}\},\qquad p\ge2,\label{eq:Kp3-main}
\end{align}
where $\cl$ means the topological closure. The parameters $(\alpha,\beta)$ are normalized as in~\eqref{eq:endpairs-main}. Together with $\Rp^3$, these make up all three-dimensional semi-homogeneous cones~\citep{HildebrandAffine2014}.

%% file: sections/section_03_the_classification.tex
\section{The classification} \label{sec:3}

In this section the main results of the paper are summarized. For convenience of the reader, the cones are collected in two tables, while a third table lists the non-trivial generators of the Lie algebra of its automorphism subgroup $G$.

\begin{theorem}[Complete classification]\label{thm:classification}
Let $K\subset\R^4$ be a semi-homogeneous proper convex cone. Then, after a linear change of coordinates $K$ is one of the cones in Tables~\ref{tab:classification-known} and~\ref{tab:classification-lorentz}. Conversely, every cone in these tables is proper and semi-homogeneous. Apart from the overlaps explicitly stated in the tables, no two displayed cones of distinct types or with distinct parameter values are linearly isomorphic.
\end{theorem}

By Lemma \ref{lem:suffCondition} a semi-homogeneous cone in $\mathbb R^4$ is either homogeneous or has principal orbits of dimension 3. The homogeneous cones are known \citep{KaneyukiTsuji74}, while the non-homogeneous ones split in two cases (Section \ref{sec:4}) which are considered in Sections \ref{sec:5} and \ref{sec:6}, respectively. 

\begin{table}[h]
\caption{Homogeneous, product, and abelian families. %The notation follows the explicit-cone conventions used in the optimization paper.
}
\label{tab:classification-known}
\centering
\small
\setlength{\tabcolsep}{2.2pt}
\renewcommand{\arraystretch}{1.16}
\begin{tabularx}{\textwidth}{@{}p{1.02in}|X| p{1.33in}@{}}
\toprule
family & cone & parameters and overlaps\\
\midrule
homogeneous & $\Rp^4$, $L_3\times\Rp$, $L_4$ & mutually nonisomorphic\\
products & $K_{\exp}^{3}\times\Rp$, $K_{\alpha,\beta}^{p,3}\times\Rp$ & $p\ge2$, $(\alpha,\beta)\in\EndPairs$; $K_{1,1}^{2,3}\times\Rp\simeq L_3\times\Rp$\\
asymmetric power &
$\displaystyle K_{\alpha,\beta}^{p,q,r}=\cl\!\left\{\begin{array}{c}(x,y,z,w)\in\R^4: x,y,z>0,\\[0mm] -\alpha x^{1/p}y^{1/q}z^{1/r}<w<\beta x^{1/p}y^{1/q}z^{1/r}\end{array}\right\}$
& $p,q,r>1$, $p^{-1}+q^{-1}+r^{-1}=1$; unordered $\{p,q,r\}$; $(\alpha,\beta)\in\EndPairs$\\
rotational-power &
$\displaystyle K_{\lambda}^{2,2}=\cl\left\{\begin{array}{c} (x,y,z,w)\in\R^4:x,y>0, \\[0mm] \sqrt{z^2+w^2}<x^\lambda y^{1-\lambda}\end{array}\right\}$
& $0<\lambda<\tfrac12$; $\lambda=\tfrac12$ gives $L_4$\\
exponential-power &
$\displaystyle K_{\mathrm{pow},a}=\cl\left\{\begin{array}{c}(x,y,z,w)\in\R^4:z,w>0, \\[0mm] y/z>\e^{x/z}(w/z)^{-a}\end{array}\right\}$
& $0<a\le1$\\
exponential-quadratic &
$\displaystyle K_{\mathrm{par}}=\cl\left\{\begin{array}{c}(x,y,z,w)\in\R^4:z>0,\\[0mm]y/z>\exp(x/z+(w/z)^2)\end{array}\right\}$
& no parameter\\
\bottomrule
\end{tabularx}
\end{table}

\begin{table}[h]
\caption{Families produced by the nonabelian matrix equation.}
\label{tab:classification-lorentz}
\centering
\small
\setlength{\tabcolsep}{2.1pt}
\renewcommand{\arraystretch}{1.17}
\begin{tabularx}{\textwidth}{@{}p{1.03in}|X| p{1.28in}@{}}
\toprule
family & cone & parameters and overlaps\\
\midrule
asymmetric mixed-power &
$\displaystyle K^N_{\gamma,\alpha,\beta}=\cl\left\{\begin{array}{c}(x,y,z,w)\in\R^4:x,z>0,\Delta>0,\\[0mm]-\alpha H_\gamma<w<\beta H_\gamma\end{array}\right\}$
& $0<\gamma<2$, $(\alpha,\beta)\in\EndPairs$; $K^N_{1,1,1}\simeq L_4$, $(\alpha,\beta) = (0,+\infty)$ is $L_3 \times \mathbb R_+$\\
cylindrical wedge & $\displaystyle K^N_{0,1,1}=\cl\left\{\begin{array}{c}(x,y,z,w)\in\R^4:x,z>0,\Delta>0,\\[0mm]-z<w<z\end{array}\right\}$ & cones $K^N_{0,\alpha,\beta}$: $\alpha,\beta < +\infty$ is isomorphic to $K^N_{0,1,1}$; $(\alpha,\beta)=(0,+\infty)$ and $(\alpha,\beta) = (1,+\infty)$ give $L_3\times\Rp$\\
croissant & $K^N_{2,1,1}=\cl\left\{\begin{array}{c}(x,y,z,w)\in\R^4:x,z>0,\Delta>0,\\[0mm]-\Delta<wz<\Delta\end{array}\right\}$ & cones $K^N_{2,\alpha,\beta}$: $0 < \alpha,\beta < +\infty$ is isomorphic to $K^N_{2,1,1}$; $(\alpha,\beta)=(0,+\infty)$ and $(\alpha,\beta)=(0,1)$ give $L_3\times\Rp$; $(\alpha,\beta) = (1,\infty)$ is $K^N_{2,1,1}$\\
one-sided mixed-power &
$\displaystyle K^N_{\mathrm{mix},\gamma}=\cl\{(x,y,z,w)\in\R^4:x,z>0,\Delta>0,w>H_\gamma\}$
& $\gamma<0$\\
exterior mixed-power &
$\displaystyle K^N_{\mathrm{ext},\gamma}=\cl\left\{\begin{array}{c}(x,y,z,w)\in\R^4:z>0,\\[0mm]w>z^{1-\gamma}(y^2-2xz)_+^{\gamma/2}\end{array}\right\}$
& $\gamma>2$\\
logarithmic-parabolic &
$\displaystyle K^N_{\log}=\cl\left\{\begin{array}{c}(x,y,z,w)\in\R^4:x,z>0,\Delta>0,\\[0mm]w<z\log(\Delta/z^2)\end{array}\right\}$
& no parameter\\
exponential-parabolic &
$\displaystyle K^N_{\log,2}=\cl\left\{\begin{array}{c}(x,y,z,w)\in\R^4:z,w>0,\\[0mm]\Delta>2wz\log(w/z)\end{array}\right\}$
& no parameter\\
\bottomrule
\end{tabularx}
\end{table}

The cone $K^{2,2}_{\lambda}$ is a generalized power cone in the sense of \citep{LinEtAl2024}. The names for the cones $K^N_{0,1,1}$ and $K^N_{2,1,1}$ are motivated by the shape of their affine compact sections, depicted in Section \ref{sec:naming} in the Appendix. Along with these, the compact sections of the exponential-quadratic and exponential-parabolic cones are depicted.

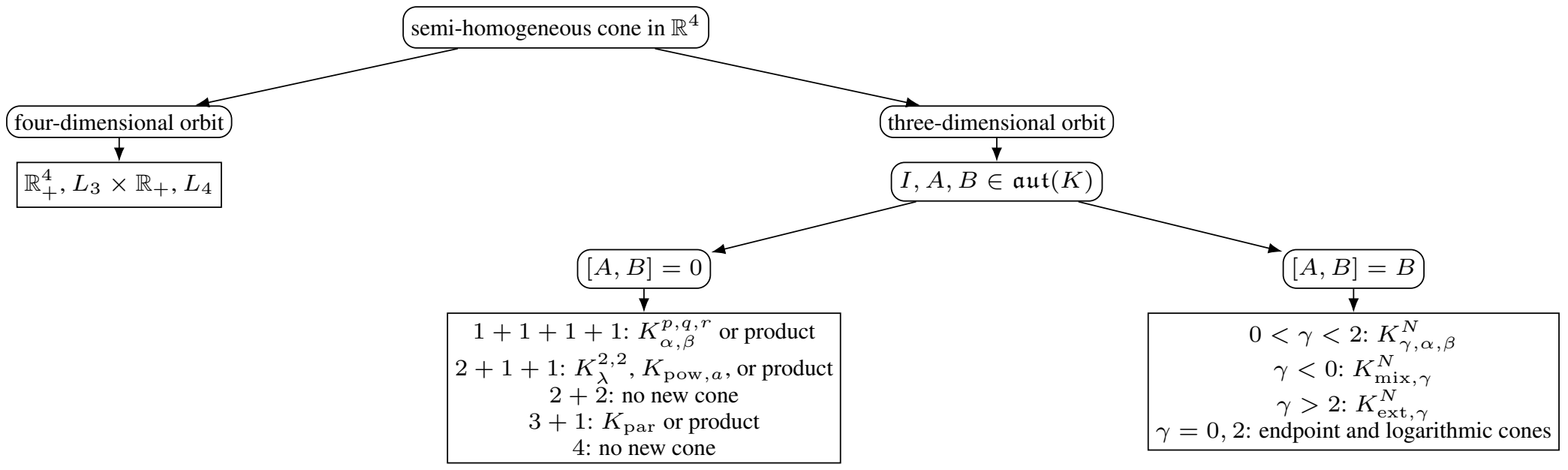
\begin{figure}[t]
\centering
\resizebox{0.98\linewidth}{!}{%
\begin{tikzpicture}[
 font=\scriptsize,node distance=2.5mm and 6mm,
 box/.style={draw,rounded corners,align=center,inner sep=2.3pt},
 leaf/.style={draw,align=center,inner sep=2.2pt},
 arr/.style={-{Latex[length=1.5mm]},thin}]
\node[box] (root) {semi-homogeneous cone in $\R^4$};
\node[box,below left=6mm and 18mm of root] (hom) {four-dimensional orbit};
\node[leaf,below=of hom] (hlist) {$\Rp^4$, $L_3\times\Rp$, $L_4$};
\node[box,below right=6mm and 18mm of root] (non) {three-dimensional orbit};
\node[box,below=of non] (red) {$I,A,B\in\aut(K)$};
\node[box,below left=5mm and 19mm of red] (ab) {$[A,B]=0$};
\node[box,below right=5mm and 19mm of red] (aff) {$[A,B]=B$};
\node[leaf,below=of ab] (ablist) {$1+1+1+1$: $K^{p,q,r}_{\alpha,\beta}$ or product\\$2+1+1$: $K^{2,2}_\lambda$, $K_{\mathrm{pow},a}$, or product\\$2+2$: no new cone\\$3+1$: $K_{\mathrm{par}}$ or product\\$4$: no new cone};
\node[leaf,below=of aff] (afflist) {$0<\gamma<2$: $K^N_{\gamma,\alpha,\beta}$\\$\gamma<0$: $K^N_{\mathrm{mix},\gamma}$\\$\gamma>2$: $K^N_{\mathrm{ext},\gamma}$\\$\gamma=0,2$: endpoint and logarithmic cones};
\draw[arr] (root)--(hom); \draw[arr] (hom)--(hlist);
\draw[arr] (root)--(non); \draw[arr] (non)--(red);
\draw[arr] (red)--(ab); \draw[arr] (red)--(aff);
\draw[arr] (ab)--(ablist); \draw[arr] (aff)--(afflist);
\end{tikzpicture}}
\caption{The complete case tree. Product branches occur when one coordinate is scaled independently or one weight is zero.}
\label{fig:tree}
\end{figure}

Table~\ref{tab:matrices-main} gives two matrices for every nonproduct family. Together with $I$, they generate a three-dimensional subgroup of automorphisms of the cone. The last column is a complete invariant of its two-dimensional action on rays. On Fig.~\ref{fig:tree} we depict the coarse classification scheme leading to the Tables \ref{tab:classification-known} and \ref{tab:classification-lorentz}.

\begin{table}[t]
\caption{Explicit matrix generators and orbit invariants.}
\label{tab:matrices-main}
\centering
\small
\setlength{\tabcolsep}{2.0pt}
\renewcommand{\arraystretch}{1.13}
\begin{tabularx}{\textwidth}{@{}p{0.90in}p{2.98in}X@{}}
\toprule
family & algebra generators $A,B$ & orbit invariant\\
\midrule
$K^{p,q,r}_{\alpha,\beta}$ &
$\displaystyle A=\begin{pmatrix}1&0&0&0\\0&0&0&0\\0&0&0&0\\0&0&0&1/p\end{pmatrix},\quad
B=\begin{pmatrix}0&0&0&0\\0&1&0&0\\0&0&0&0\\0&0&0&1/q\end{pmatrix}$
& $\displaystyle w/(x^{1/p}y^{1/q}z^{1/r})$\\
$K^{2,2}_\lambda$ &
$\displaystyle A=\begin{pmatrix}1&0&0&0\\0&0&0&0\\0&0&\lambda&0\\0&0&0&\lambda\end{pmatrix},\quad
B=\begin{pmatrix}0&0&0&0\\0&0&0&0\\0&0&0&-1\\0&0&1&0\end{pmatrix}$
& $\displaystyle \sqrt{z^2+w^2}/(x^\lambda y^{1-\lambda})$\\
$K_{\mathrm{pow},a}$ &
$\displaystyle A=\begin{pmatrix}0&0&1&0\\0&1&0&0\\0&0&0&0\\0&0&0&0\end{pmatrix},\quad
B=\begin{pmatrix}0&0&0&0\\0&-a&0&0\\0&0&0&0\\0&0&0&1\end{pmatrix}$
& $\displaystyle \log(y/z)-x/z+a\log(w/z)$\\
$K_{\mathrm{par}}$ &
$\displaystyle A=\begin{pmatrix}0&0&1&0\\0&1&0&0\\0&0&0&0\\0&0&0&0\end{pmatrix},\quad
B=\begin{pmatrix}0&0&0&-2\\0&0&0&0\\0&0&0&0\\0&0&1&0\end{pmatrix}$
& $\displaystyle \log(y/z)-x/z-(w/z)^2$\\
$K^N_{\gamma,\alpha,\beta}$, $K^N_{\mathrm{mix},\gamma}$ &
$\displaystyle A_\gamma=\begin{pmatrix}2&0&0&0\\0&1&0&0\\0&0&0&0\\0&0&0&\gamma\end{pmatrix},\quad
B_L=\begin{pmatrix}0&1&0&0\\0&0&1&0\\0&0&0&0\\0&0&0&0\end{pmatrix}$
& $\displaystyle w/(\Delta^{\gamma/2}z^{1-\gamma}),\quad \Delta>0$\\
$K^N_{\mathrm{ext},\gamma}$ &
$\displaystyle A_\gamma=\begin{pmatrix}2&0&0&0\\0&1&0&0\\0&0&0&0\\0&0&0&\gamma\end{pmatrix},\quad B_L=\begin{pmatrix}0&1&0&0\\0&0&1&0\\0&0&0&0\\0&0&0&0\end{pmatrix}$
& $\displaystyle w/[z^{1-\gamma}\Delta^{\gamma/2}],\quad \Delta>0$\\
$K^N_{\log}$ &
$\displaystyle A_{\log}=\begin{pmatrix}2&0&0&0\\0&1&0&0\\0&0&0&0\\0&0&2&0\end{pmatrix},\quad B_L=\begin{pmatrix}0&1&0&0\\0&0&1&0\\0&0&0&0\\0&0&0&0\end{pmatrix}$
& $\displaystyle w/z-\log(\Delta/z^2)$\\
$K^N_{\log,2}$ &
$\displaystyle A_{\log,2}=\begin{pmatrix}2&0&0&2\\0&1&0&0\\0&0&0&0\\0&0&0&2\end{pmatrix},\quad B_L=\begin{pmatrix}0&1&0&0\\0&0&1&0\\0&0&0&0\\0&0&0&0\end{pmatrix}$
& $\displaystyle \Delta/(2wz)-\log(w/z)$\\
\bottomrule
\end{tabularx}
\end{table}

%% file: sections/section_04_reduction_to_two_matrix_equations.tex
\section{Classification of Lie algebras} \label{sec:4}

In this section we show that in the case the principal orbits of the automorphism group have dimension 3, the stabilizer subgroup of interior points in principal orbits can be removed, in the sense that there exists a subgroup $G \subset \Aut(K)$ with trivial generic stabilizer that has the same principal orbit dimension as $\Aut(K)$.

\begin{theorem}[Three-dimensional subgroup]\label{thm:reduction-main}
Let $K\subset\R^4$ be a proper convex cone and suppose that a principal $\Aut(K)$-orbit has dimension three. Then $\Aut(K)$ has a connected three-dimensional subgroup $G$ with Lie algebra $\mathfrak g$ that contains the positive homotheties and has three-dimensional principal orbits. There exist matrices $A,B \in \mathfrak g$ such that
\[
 \mathfrak g=\Span(I,A,B),\qquad [A,B]=0\quad\text{or}\quad[A,B]=B.
\]
\end{theorem}

\begin{proof}
Choose an interior point $\hat x \in K^{\circ}$ inside a principal orbit and define the coordinate system in $\mathbb R^4$ as in Lemma \ref{lem:structureAut} in Section \ref{sec:compactStabilizers} in the Appendix. Note that the stabilizer of $\hat x$ can have dimension at most one. 

A trivial stabilizer gives the required three-dimensional subgroup $G$ immediately as the connected component of the identity in $\Aut(K)$. 

If the stabilizer is non-trivial, then Lemma \ref{lem:structureStabilizer} in Section \ref{sec:trivialStabilizer} in the Appendix applies. It gives a normalized basis $\{A,B,C\}$ of the subalgebra $\mathfrak g_0 \subset \aut(K)$ of trace-less matrices, from which we can infer its Bianchi type \cite{Bianchi}. 

The basis matrices $A,B$ contain two parameters $a_{24},a_{42}$, where $a_{ij}$ denotes the elements of $A$. Set $\kappa = -(1+a_{24}a_{42})$. The structure of $A,B,C$ implies that $\mathfrak g_0$ is isomorphic to one of the three Lie algebras $\mathfrak e(2)$ (for $\kappa = 0$), $\mathfrak{sl}_2(\R)$ (for $\kappa < 0$), and $\mathfrak{so}_3$ (for $\kappa > 0$). The first two contain a two-dimensional subgroup with a 2-dimensional orbit, which together with the homotheties yield the sought subgroup $G$.

In the last case the automorphism group $\Aut(K)$ contains the subgroup of unimodular rotations in $\mathbb R^4$ which leave the vector $(-a_{24},0,0,1)^T$ fixed. Then we must have $K\simeq L_4$, and $K$ is homogeneous, contradicting the assumptions of the theorem. 

Finally, every two-dimensional real Lie algebra is abelian or has a basis with $[A,B]=B$ \citep{Mubarakzyanov1963}.
\end{proof}

% A second elementary restriction (Lemma \ref{lem:nilpotent-obstruction}) will be used repeatedly. If a nonzero nilpotent matrix $N\in\aut(K)$ has an even-dimensional largest Jordan block, then for a generic $v\in K^\circ$ the highest term in $\e^{tN}v$ has odd degree in $t$. Positive normalization as $t\to+\infty$ and $t\to-\infty$ places two mutually opposite nonzero vectors in $K$, which contradicts pointed-ness. In dimension four this excludes every nonzero matrix $A$ satisfying $A^2 = 0$ (a square-zero matrix) and the matrix with ones in positions $(1,2),(2,3),(3,4)$ (a 4-dimensional nilpotent Jordan block) and similar matrices.

% On a connected open set where the two-dimensional action on rays has a scalar invariant $q$ and the level sets of $q$ are its orbits, every connected invariant interior is $q^{-1}(J)$ for an interval $J$. This follows simply because the image of a connected set is an interval and invariance contains every full level set. Once $J$ is known, convexity is checked directly from the defining inequalities.

%% file: sections/section_05_the_abelian_case.tex
\section{The abelian case} \label{sec:5}

In this section we consider the case when the orbits of $\Aut(K)$ are generated by a 3-dimensional abelian subgroup $G \subset \Aut(K)$ whose Lie algebra is spanned by commuting matrices $\{I,A,B\}$.

\begin{theorem}[Simultaneous invariant subspace decomposition]\label{thm:weights-main}
Let $K\subset\R^4$ be proper, let $A,B$ commute, and assume
\[
 \Span(I,A,B)\subset\aut(K).
\]
The space $\R^4$ is a direct sum of common real invariant subspaces for $A,B$ with the full dimension decomposing by one of the patterns
\begin{equation}\label{eq:partitions-main}
 1+1+1+1,
 \qquad2+1+1,
 \qquad2+2,
 \qquad3+1,
 \qquad4.
\end{equation}

Suppose the principal orbits of the group generated by $I,A,B$ are 3-dimensional. Then on a two-dimensional common subspace, after subtracting scalar matrices (multiples of the identity), the restrictions of $A,B$ can be made multiples of one of
\[
 \begin{pmatrix}0&1\\0&0\end{pmatrix},
 \qquad
 \begin{pmatrix}0&-1\\1&0\end{pmatrix}.
\]
On a three-dimensional common subspace, after choosing a nonzero nilpotent generator, the two restrictions can be made polynomials of degree at most two in
\[
 J_3 = \begin{pmatrix}0&1&0\\0&0&1\\0&0&0\end{pmatrix}.
\]
\end{theorem}

The proof is provided in Section \ref{sec:ProofTheorem3} in the Appendix. It relies on the canonical form of pairs of matrices with respect to similarity transformations. Next we consider the emerged subcases one by one. Technical details are delegated to the remaining subsections of Section \ref{app:abelian} in the Appendix. %The four-dimensional case is treated in Appendix~\ref{app:abelian}.

\subsection{Four one-dimensional common eigenspaces}\label{sec:5-1}

There is one common affine relation $\sum_{k=1}^4 \alpha_j\lambda_j = 0$, $\sum_j \alpha_j = 0$ on the four pairs of eigenvalues $\lambda^A_j,\lambda^B_j$ of $A,B$. If one coefficient $\alpha_j$ has sign opposite to the other three, a change of generators gives the basis
\[
A=\begin{pmatrix}1&0&0&0\\0&0&0&0\\0&0&0&0\\0&0&0&1/p\end{pmatrix},\quad B=\begin{pmatrix}0&0&0&0\\0&1&0&0\\0&0&0&0\\0&0&0&1/q\end{pmatrix},
\quad
I-A-B = \begin{pmatrix}0&0&0&0\\0&0&0&0\\0&0&1&0\\0&0&0&1/r\end{pmatrix},
\]
where $p^{-1}+q^{-1}+r^{-1}=1$. The orbit invariant is given by
\[
 f=\frac{w}{|x|^{1/p}|y|^{1/q}|z|^{1/r}}.
\]
In every one of the 16 open orthants of $\mathbb R^4$, the orbits of the group action are exactly the level sets of $f$. The coordinates $x,y,z$ cannot assume opposite signs in $K$, because otherwise the corresponding canonical basis vectors are in $K$ together with their negation. Hence we may assume without loss of generality that $x,y,z \geq 0$ on $K$. Then the interior $K^{\circ}$ consists of a union of orbits on which the value of $f$ lies in an open interval. 

Let
\[
 m(x,y,z)=x^{1/p}y^{1/q}z^{1/r}.
\]
For $h=(h_1,h_2,h_3)$,
\begin{equation}\label{eq:geom-mean-hessian}
 D^2m[h,h]
 =m\left[\left(\sum_i\lambda_i\frac{h_i}{x_i}\right)^2
 -\sum_i\lambda_i\left(\frac{h_i}{x_i}\right)^2\right]\le0
\end{equation}
by weighted Cauchy--Schwarz. Thus $m$ is concave. The inequality $w\le\beta m$ is convex exactly for $\beta\ge0$, and $w\ge-\alpha m$ is convex exactly for $\alpha\ge0$. Reflection and positive scaling of $w$ reduce all nontrivial intervals to~\eqref{eq:endpairs-main}. This gives the cones $K^{p,q,r}_{\alpha,\beta}$.

In the case of two positive and two negative coefficients $\alpha_j$ no invariant cone can exist. If the affine relation has zero coefficients, then the invariant cones must be direct products. A detailed analysis of  these two cases is given in Appendix~\ref{sec:B2}.

\subsection{One two-dimensional and two one-dimensional common subspaces}

Let the entire space be decomposed into a sum $\mathbb{R}^4 = V \oplus \mathbb{R} \oplus \mathbb{R}$, and let the generators be written in the corresponding basis.

If both generators $A, B$ act on a two dimensional subspace as $I_2$, then this can be attributed to the previous subsection.

The case where $A,B$ generate action on $V$ as a shear and a rotation is impossible because then they do not commute.

% We can assume that $A_{11} \neq A_{22}$, otherwise $diag(0,0,1,1) \in \mathfrak{g}$ and by Lemma \ref{lem:product2D} the cone is represented as a product.
% добавлю в аппендикс

If an element of the algebra generates a rotation on $V$, then $A, B$ can be chosen of the form \eqref{eq:norm-matrices-app}. A detailed explanation of this is contained in Section~\ref{sec:abelian211rotation} of Appendix. Lemma
\ref{lem:normPowerCone} shows that in this case only $K^{2,2}_\lambda$ appears.

When the generated action on $V$ is a shear, $A, B$ can be normalized to the form \eqref{eq:Kpow-matrices-app}. The proof is provided in Section \ref{sec:B3.2}. Lemma \ref{lem:K_pow_mu} proves that only
$K_{\mathrm{pow},a}$ appears in this case.

\subsection{Two two-dimensional common subspaces}

The entire space is the sum of invariant subspaces $\mathbb{R}^4 = V \oplus W$. Let the generators be written in the corresponding basis.

If $A,B$ act:
\begin{itemize}
    \item both  as $I_2$ on $V$ or $W$, then this can be attributed to the previous sections.
    \item as a shear on $V$, scalar on $W$ for $A$ and a shear on $W$, scalar on $V$ for $B$, then algebra contains two elements:
    \[
     A=\begin{pmatrix}a&1&0&0\\0&a&0&0\\0&0&0& 0\\0&0&0&0\end{pmatrix},\quad
     B=\begin{pmatrix}0&0&0&0\\0&0&0&0\\0&0&b&1\\0&0&0&b\end{pmatrix},
    \]
    Case $a = 0$ or $b = 0$ is impossible by lemma \ref{lem:nilpotent-obstruction} . In the remaining case, as in lemma \ref{lem:product2D}, we obtain that $K$ is a product cone.
    \item both as a shear on $V$ and $W$, then we do not fall into the previous case only if pairs of shear coefficients of $A$ and $B$ are dependent, but at that rate algebra contains $\diag\{1, 1, 0, 0\}$ and by lemma \ref{lem:product2D} $K$ is a product.
    \item as a rotation on $V$, scalar on $W$, then relative scaling isolates a full circle and the resulting cone cannot be pointed.
\end{itemize}

The two remaining mixed pairs are
\[
 A=\begin{pmatrix}a&-1&0&0\\1&a&0&0\\0&0&0&0\\0&0&0&0\end{pmatrix},\quad
 B=\begin{pmatrix}b&0&0&0\\0&b&0&0\\0&0&0&1\\0&0&0&0\end{pmatrix},
\]
and
\[
 A=\begin{pmatrix}0&-1&0&0\\1&0&0&0\\0&0&0&1\\0&0&0&0\end{pmatrix},\quad
 B=\begin{pmatrix}1&0&0&0\\0&1&0&0\\0&0&0&0\\0&0&0&0\end{pmatrix}.
\]
No new cones appear in this case. A detailed argument and analysis of other cases is given in Section \ref{subsec:mixed22} of Appendix.

\subsection{A three-dimensional common subspace and one line}

After removing scalar terms and a change of generators, the matrices are
\[
 A=\begin{pmatrix}0&1&0&0\\0&0&1&0\\0&0&0&0\\0&0&0&0\end{pmatrix},
 \qquad
 B=\begin{pmatrix}0&0&c&0\\0&0&0&0\\0&0&0&0\\0&0&0&1\end{pmatrix}.
\]
The orbit through $(0,0,1,1)$ is
\[
 \{(s^2/2+c\log w,s,1,w):s\in\R,\ w>0\}.
\]
Thus in the section $z=1$ the boundary graph is $x=y^2/2+c\log w$. For $c>0$ its Hessian has one positive and one negative eigenvalue, so no convex side exists. For $c=0$ the cone splits and is isomorphic to $L_3\times\Rp$. For $c<0$, rescale to $c=-1/2$ and apply
\[
 (x,y,z,w)\longmapsto(-2x,w,z,y).
\]
The graph becomes $y=\exp(x+w^2)$, and the resulting cone is $K_{\mathrm{par}}$. A full argument is given in Section \ref{sec:B5} of Appendix.

\subsection{One four-dimensional common subspace}

The only real case surviving the nilpotent obstruction is generated by
\[
 A=\begin{pmatrix}0&1&0&0\\0&0&0&1\\0&0&0&0\\0&0&0&0\end{pmatrix},
 \qquad
 B=\begin{pmatrix}0&0&1&0\\0&0&0&0\\0&0&0&1\\0&0&0&0\end{pmatrix}.
\]
It preserves $2xw-y^2-z^2$. The orbit in $w=1$ is
\[
 \{((s^2+t^2+c)/2,s,t,1):s,t\in\R\}.
\]
Only the epigraph side is convex; translating $x$ normalizes $c=0$ and gives $L_4$. The nonreal four-dimensional case contains a nonzero square-zero matrix or has only a one-dimensional action on rays. A full argument is given in Section \ref{sec:B6} of Appendix. This completes the abelian list.

%% file: sections/section_06_the_nonabelian_case.tex
\section{The nonabelian case} \label{sec:6}

\begin{proposition}[Complete solution of the affine matrix equation]\label{prop:affine-main}
Assume $[A,B]=B$ and $\Span(I,A,B)\subset\aut(K)$ for a proper cone. After a real change of coordinates and a change of the two projective generators, we may assume
\begin{equation}\label{eq:affine-main}
 B=\begin{pmatrix}0&1&0&0\\0&0&1&0\\0&0&0&0\\0&0&0&0\end{pmatrix},
 \qquad
 A=\begin{pmatrix}2&0&c&d\\0&1&0&0\\0&0&0&0\\0&0&e&\gamma\end{pmatrix}.
\end{equation}
If $\gamma\ne0,2$, then all three entries $c,d,e$ can be set to 0. If $\gamma=0$, the entry $e$ may remain. If $\gamma=2$, the entry $d$ may remain.
\end{proposition}

Appendix~\ref{sec:commutatorSystem} contains the full derivation of the proposition. The core part of the proof is a direct solution of the matrix equation $[A,B] = B$ in Lemma \ref{lem:commutator-solution}. First a finite list of canonical forms for $B$ is found, then the linear system on $A$ is solved. 

The proposition yields a number of cases, each with a parameter. In Sections \ref{sec:C2} and \ref{sec:C3} in the appendix each of the cases is considered. The invariant convex cones which arise are classified. These are the cones presented in Table \ref{tab:classification-lorentz}. This finishes the nonabelian list.

%\subsection{Why the list is complete at the level of cones}

%Appendices~\ref{app:abelian} and~\ref{app:affine} write every principal orbit and its separating invariant explicitly. Lemma~\ref{lem:regular-boundary-direct} supplies a regular boundary orbit; normalized-flow limits fix the abelian sign chambers; and Lemma~\ref{lem:affine-chambers} shows that crossing $\Delta=0$ yields only $K^N_{\mathrm{ext},\delta}$ or a product. The invariant-interval lemma and the displayed convexity or midpoint calculations therefore exhaust all interiors and all finite endpoints. Closedness forbids adding independent boundary pieces, while the product and nilpotent lemmas remove the remaining cases. Appendix~\ref{app:verification} computes the full algebras and normalizers, including $\gamma=1,-2,4$, and Appendix~\ref{app:completion} assembles the complete proof.

%% file: sections/section_07_conclusion.tex
\section{Conclusion} \label{sec:8}

In this work we completely classified the 4-dimensional semi-homogeneous cones up to linear isomorphism. Semi-homogeneous cones are defined by having a sufficiently rich automorphism group, such that its principal orbits are of codimension at most one. 

Semi-homogeneous cones are prospective candidates for building blocks in conic optimization, because the computation of the canonical barrier on them depends on the solution merely of an ordinary differential equation rather than the partial differential Monge-Amp\`ere equation which originally defines the barrier. The derivation and solution of these ODEs and the construction of the canonical barriers will be the subject of a separate paper.

In the resulting classification there appear families of cones which are already known, such as the generalized power cones, but also new families of cones and new singleton cones not belonging to a continuous family.

%% file: sections/appendix_structure.tex
\section{Structural lemmas}\label{app:structure}

\subsection{The canonical potential and compact stabilizers} \label{sec:compactStabilizers}

\begin{lemma}[Normalization at an interior point]\label{lem:structureAut}
    Let $K \subset \mathbb R^4$ be a proper convex cone with canonical barrier $F$. Let $\hat x$ be an interior point in an orbit of $\Aut(K)$ of dimension 3. Let $\mathfrak g_0$ be the trace-less subalgebra of the Lie algebra $\mathfrak g$ of $\Aut(K)$. Then there exists a linear coordinate system in $\mathbb R^4$ such that 
    \begin{itemize}
        \item $\hat x = e_1$
        \item the orbit of $\hat x$ is orthogonal to $e_4$ at $\hat x$
        \item every element of $\mathfrak g_0$ has the form $\begin{pmatrix}
            0 & a_{12} & a_{13} & 0 \\ a_{12} & a_{22} & a_{23} & a_{24} \\ a_{13} & a_{32} & a_{33} & a_{34} \\ 0 & a_{42} & a_{43} & a_{44}
        \end{pmatrix}$
        \item every element of the Lie algebra of the stabilizer subgroup of $\hat x$ is proportional to $\begin{pmatrix}
            0 & 0 & 0 & 0 \\ 0 & 0 & -1 & 0 \\ 0 & 1 & 0 & 0 \\ 0 & 0 & 0 & 0
        \end{pmatrix}$
    \end{itemize}
    Here $e_i$ are the basis vectors of the coordinate system in $\mathbb R^4$.
\end{lemma}

\begin{proof}
Choose the coordinate system in $\mathbb R^4$ such that  $\hat x = e_1$ and such that the Hessian of the canonical potential is given by $F''(\hat x) = aI_4$ ($I_4$ being the identity matrix) for some $a > 0$. Since for logarithmically homogeneous functions we have $F''(x)x = -F'(x)$ \citep{NesterovNemirovskii1994}, the level surface of $F$ is orthogonal to $e_1$ at $\hat x$. Rotate the coordinate system such that the orbit of $\hat x$ becomes orthogonal to $e_4$ at $\hat x$. Then the tangent space to the intersection of the orbit with the level surface at $\hat x$ is spanned by $e_2,e_3$.

Let now $A \in \mathfrak g_0$ and $g_t = e^{tA}$. Then $\det\,g_t = 1$ and $g_t$ preserves the level surfaces of $F$ for all $t$ by~\eqref{automorphismSymmetry}. Differentiating~\eqref{automorphismSymmetry} with respect to $x$ at $x = \hat x$, we obtain
\[ e^{tA^T}F'(e^{tA}\hat x) = F'(\hat x)\quad \forall\ t.
\]
For small $t$ we have
\[ F'(e^{tA}\hat x) = F'(\hat x) + F''(\hat x)(e^{tA}-I)\hat x + O(t^2) = -F''(\hat x)\hat x + tF''(\hat x)A\hat x + O(t^2),
\]
and hence
\[ (I + tA^T + O(t^2))(F''(\hat x)(-I + tA)\hat x + O(t^2)) = -F''(\hat x)\hat x.
\]
It finally follows by comparing the coefficients at $t$ that
\begin{equation} \label{symmetricCondition}
(A^T - A)\hat x = 0.
\end{equation} 

On the other hand, since $g_t\hat x = \hat x + tA\hat x + O(t^2)$ has to be in the orbit and in the level surface of $\hat x$ for all $t$, we get that $A\hat x$ is in the span of $e_2,e_3$. Together with \eqref{symmetricCondition} this proves the claim on the structure of $A$ in the lemma.

Let now $g \in \Aut(K)$ be in the connected component of the identity in the stabilizer subgroup of $\hat x$. Differentiating~\eqref{automorphismSymmetry} twice with respect to $x$ at $x = \hat x$, we get that
\[ g^TF''(\hat x)g = F''(\hat x),
\]
i.e., the automorphism $g$ preserves $F''(\hat x) = aI_4$ as a quadratic form.

Hence the connected component of the identity in the stabilizer of $\hat x$ is a subgroup of $SO(4)$. Since $g\hat x = \hat x$, the action of $g$ must amount to a 2D rotation. By~\eqref{automorphismSymmetry} $g$ preserves as a set the intersection of the orbit with the level set of $F$ through $\hat x$, which is two-dimensional. Hence $g$ must be a 2D rotation in the plane tangent to this intersection, i.e., spanned by $e_2,e_3$. The last claim of the lemma now easily follows.
\end{proof}

\subsection{Presence of a non-trivial stabilizer} \label{sec:trivialStabilizer}

\begin{lemma}[Structure of the Lie algebra when the stabilizer is non-trivial]\label{lem:structureStabilizer}
    Let $K \subset \mathbb R^4$ be a proper convex cone and suppose $\Aut(K)$ has dimension $>3$ and has an orbit of dimension 3 in the interior of $K$. Then the coordinate system in $\mathbb R^4$ can be chosen such that $e_1 \in K^{\circ}$, and the Lie algebra $\aut(K)$ of $\Aut(K)$ has dimension 4 and is spanned by the basis matrices
    \begin{equation} \label{eq:stabilizerBasis}
        I_4,\ A = \begin{pmatrix}
    0 & 1 & 0 & 0 \\ 1 & 0 & 0 & a_{24} \\ 0 & 0 & 0 & 0 \\ 0 & a_{42} & 0 & 0
\end{pmatrix}, \ B = \begin{pmatrix}
    0 & 0 & 1 & 0 \\ 0 & 0 & 0 & 0 \\ 1 & 0 & 0 & a_{24} \\ 0 & 0 & a_{42} & 0
\end{pmatrix}, \quad C = \begin{pmatrix}
            0 & 0 & 0 & 0 \\ 0 & 0 & -1 & 0 \\ 0 & 1 & 0 & 0 \\ 0 & 0 & 0 & 0
        \end{pmatrix}.
    \end{equation}
    Here $a_{24},a_{42}$ are parameters. 
    %$\kappa = -(1+a_{24}a_{42})$
\end{lemma}

\begin{proof}
Choose an interior point $\hat x \in K^{\circ}$ inside an orbit of dimension 3 and define the coordinate system in $\mathbb R^4$ as in Lemma \ref{lem:structureAut}. Note that the stabilizer of $\hat x = e_1 \in K^{\circ}$ can have dimension at most one by this lemma. 

Since $\dim\,\Aut(K) > 3$, the stabilizer of $\hat x$ must be one-dimensional. This implies that the subalgebra $\mathfrak g_0 \subset \aut(K)$ of trace-less matrices is 3-dimensional. We shall now deduce a normalized basis $\{A,B,C\}$ of $\mathfrak g_0$. Since the stabilizer is non-trivial, the Lie algebra $\mathfrak g_0$ contains an element $C$ as in the formulation, by Lemma \ref{lem:structureAut}. We choose the other two basis elements $A,B$ orthogonal to $C$, i.e., we assume $a_{23} = a_{32}$, $b_{23} = b_{32}$. Since the tangent space to the intersection of the orbit with the level surface of $\hat x$ is spanned by $e_2,e_3$, we may choose $A,B$ such that $A$ moves $\hat x$ in the direction of $e_2$, while $B$ moves it in the direction of $e_3$. Thus
\[ A = \begin{pmatrix}
    0 & 1 & 0 & 0 \\ 1 & a_{22} & a_{23} & a_{24} \\ 0 & a_{23} & a_{33} & a_{34} \\ 0 & a_{42} & a_{43} & a_{44}
\end{pmatrix}, \quad B = \begin{pmatrix}
    0 & 0 & 1 & 0 \\ 0 & b_{22} & b_{23} & b_{24} \\ 1 & b_{23} & b_{33} & b_{34} \\ 0 & b_{42} & b_{43} & b_{44}
\end{pmatrix}
\]
by virtue of Lemma \ref{lem:structureAut}.

Computing the commutators of $A,B,C$, we find that $[A,B]_{12} = -[A,B]_{21}$, $[A,B]_{13} = -[A,B]_{31}$. But then by Lemma~\ref{lem:structureAut} we must have $[A,B]_{12} = [A,B]_{13} = 0$, and it follows that $[A,B] = \kappa C$ for some $\kappa \in \mathbb R$. Further $[A,C]_{21} = [A,C]_{12} = 0$, $[A,C]_{31} = [A,C]_{13} = -1$, $[A,C]_{23} = [A,C]_{32}$ implies $[A,C] = -B$. In a similar way we get $[B,C] = A$. Inserting the above values of $A,B,C$, we obtain the expressions for $A,B$ given in the formulation of the lemma.
\end{proof}

\subsection{Homogeneous cones}

In this section we show that in the non-homogeneous case, the principal orbits have dimension 3.

\begin{lemma}[Sufficient condition for homogeneity]
\label{lem:suffCondition}
Let $K \subset \mathbb{R}^4$ be a proper convex cone whose principal
$\Aut(K)$-orbits have dimension $4$. Then $K$ is homogeneous.
\end{lemma}

\begin{proof}
Let $K^\circ$ denote the interior of $K$ and let $G=\Aut(K)^\circ$
be the identity component of the automorphism group. Since passing to the
identity component does not change orbit dimensions, $G$ has a
four-dimensional orbit $O=G\cdot x\subset K^\circ$
for some $x\in K^\circ$.

First, $O$ is open in $K^\circ$. Indeed, the differential of the orbit map
\[
\Phi_x:G\to K^\circ,
\qquad
\Phi_x(g)=gx,
\]
at the identity has image
\[
\operatorname{im} d(\Phi_x)_e=T_xO.
\]
Since $\dim O=4=\dim K^\circ$, the differential is surjective. By
$G$-equivariance it is surjective at every point of $G$, and hence
$\Phi_x$ is a submersion. Therefore $O$ is open in $K^\circ$.

It remains to show that $O$ is closed in $K^\circ$. Let $F$ be the
canonical potential of $K$ and put
\[
H(y)=D^2F(y), \qquad y\in K^\circ.
\]
For every $g\in G$, the transformation law
\[
F(gy)=F(y)-\log |\det g|
\]
implies, after differentiating twice,
\[
g^{T}H(gy)g=H(y).
\tag{*}
\]

Suppose that $g_nx\longrightarrow y\in K^\circ$
for a sequence $g_n\in G$. By continuity of $H$ and positive
definiteness of $H(y)$, there exist constants $m,M>0$ such that, for all
sufficiently large $n$,
\[
mI\preceq H(g_nx)\preceq MI.
\]
Likewise, since $H(x)\succ0$, there exist $m_0,M_0>0$ such that
\[
m_0I\preceq H(x)\preceq M_0I.
\]

Applying $(*)$ with $g=g_n$ and $y=x$, for every $v\in\mathbb{R}^4$ we
obtain
\[
v^TH(x)v
=
(g_nv)^TH(g_nx)(g_nv).
\]
Consequently,
\[
m\|g_nv\|^2
\le
v^TH(x)v
\le
M_0\|v\|^2,
\]
and therefore
\[
\|g_nv\|
\le
\sqrt{\frac{M_0}{m}}\,\|v\|.
\]
Similarly,
\[
m_0\|v\|^2
\le
v^TH(x)v
\le
M\|g_nv\|^2,
\]
so that
\[
\|g_nv\|
\ge
\sqrt{\frac{m_0}{M}}\,\|v\|.
\]
Hence the sequence $(g_n)$ is uniformly bounded, while its smallest
singular values are uniformly bounded away from zero. Passing to a
subsequence, we obtain $g_n\longrightarrow g$ for some $g\in\mathrm{GL}_4(\R)$.

The group $\Aut(K)$ is closed in $\mathrm{GL}_4(\R)$, and its identity
component $G$ is closed as well. Hence $g\in G$. Therefore
\[
y=\lim_{n\to\infty}g_nx=gx\in O.
\]
Thus $O$ is closed in $K^\circ$.

Finally, $K^\circ$ is convex and hence connected. Since $O$ is a
nonempty subset of $K^\circ$ which is both open and closed, we conclude
that $O=K^\circ$.
Thus $G$, and therefore $\Aut(K)$, acts transitively on $K^\circ$.
Hence $K$ is homogeneous.
\end{proof}

As a consequence, no non-homogeneous cone can have an orbit of dimension 4.

\subsection{Three elementary tools}

\begin{lemma}[Nilpotent obstruction]\label{lem:nilpotent-obstruction}
Let $0\ne N\in\aut(K)$ be nilpotent. Its largest Jordan block cannot have even dimension.
\end{lemma}

\begin{proof}
For a generic $v\in K^\circ$, let $d$ be maximal with $N^dv\ne0$. If the largest block has even dimension $2k$, then $d = 2k-1$ is odd. Since $\e^{tN}v\in K^{\circ}$ for every real $t$,
\[
 |t|^{-d}\e^{tN}v\longrightarrow \frac{N^dv}{d!}
 \quad(t\to+\infty),
 \qquad
 |t|^{-d}\e^{tN}v\longrightarrow-\frac{N^dv}{d!}
 \quad(t\to-\infty).
\]
Closedness places both opposite nonzero vectors in $K$, contradicting pointedness.
\end{proof}

In dimension four, Lemma~\ref{lem:nilpotent-obstruction} excludes every nonzero square-zero matrix ($N \not= 0$, $N^2 = 0$) and
\[
J_4 = \begin{pmatrix}0&1&0&0\\0&0&1&0\\0&0&0&1\\0&0&0&0\end{pmatrix}
\]
up to similarity. The only nonzero nilpotent form that up to similarity can occur as the matrix $B$ in $[A,B]=B$ is
\begin{equation}   \label{eq:three-matrix-app}
\begin{pmatrix}0&1&0&0\\0&0&1&0\\0&0&0&0\\0&0&0&0\end{pmatrix}.
\end{equation} 
Moreover, the commuting pair generated by this matrix and its square is impossible because the square is nonzero and square-zero.

\begin{lemma}[Convex graph criterion]\label{lem:graph-convexity}
Let a $C^2$ surface in an affine section be written locally as $x=f(y,z)$. If it is a boundary patch of a convex set lying on the side $x\ge f(y,z)$, then $f$ is convex. If the convex set lies on the side $x\le f(y,z)$, then $f$ is concave. In particular, if $D^2f$ has one positive and one negative eigenvalue at a point, neither side is locally convex.
\end{lemma}

\begin{proof}
Intersect the set with every affine plane obtained by fixing a line in the $(y,z)$-plane. The resulting epigraph or hypograph in that plane is convex, which is equivalent to convexity or concavity of the restriction of $f$ to every line. This is equivalent to the asserted sign of $D^2f$.
\end{proof}

\subsection{A regular boundary orbit exists}

\begin{lemma}[Boundary point of full projective rank]\label{lem:regular-boundary-direct}
For every effective matrix pair used in Appendices~\ref{app:abelian} and~\ref{app:affine}, a full-dimensional invariant proper cone is either one of the product or homogeneous cases already listed, or its boundary contains a point $v$ for which
\[
 \rank(v,Av,Bv)=3.
\]
Consequently, a relatively open part of the boundary is a two-dimensional orbit on rays.
\end{lemma}

\begin{proof}
Projectivize $K^\circ$ to a bounded convex domain in an affine chart. Its boundary has dimension two. Projective rank drops exactly where all $3\times3$ minors of
\[
 M(v)=(v,Av,Bv)
\]
vanish.

For the mixed $2+2$ matrices in~\eqref{eq:mixed-independent-matrices-app} and~\eqref{eq:mixed-proportional-matrices-app}, the four maximal minors are, respectively,
\[
 w^2(ay+x),\quad w^2(ax-y),\quad-bw(x^2+y^2),\quad-(x^2+y^2)(bz-w),
\]
and
\[
 -w^2y,\quad-w^2x,\quad-w(x^2+y^2),\quad-z(x^2+y^2).
\]
Their common zero sets are finite unions of projective curves. The same direct multiplication gives only projective curves for the shear case and the $3+1$ case. For the commuting four-dimensional matrices~\eqref{eq:four-final-matrices-app}, the only two-dimensional component of the rank-drop set is the plane $w=0$. A bounded three-dimensional convex body cannot have its entire boundary in one affine plane.

For $A_\gamma,B_L$, the four maximal minors, up to nonzero common factors, are
\begin{equation}\label{eq:affine-rank-minors}
 \gamma wz^2,
 \qquad \gamma wyz,
 \qquad-w\bigl((\gamma-2)xz+(1-\gamma)y^2\bigr),
 \qquad z\Delta.
\end{equation}
For $\gamma\notin\{0,1,2\}$ their common zero set has projective dimension at most one. At $\gamma=1$ the only additional component is the plane $z=0$, which cannot contain the whole boundary. At $\gamma=2$ the additional pieces are intersections of planes and quadrics and again have dimension at most one. At $\gamma=0$ the only possible surface component is $\Delta=0$, independent of $w$. If the whole boundary were contained there, every projective segment obtained by varying $w$ over a fixed null base point would remain in the boundary; its conic lift would contain a complete line, contrary to pointedness. For the two logarithmic matrices, substitution in $(v,Av,Bv)$ gives only projective curves in the rank-drop set.

For four common eigenvectors, a boundary contained in the union of the coordinate hyperplanes gives a polyhedral diagonal-invariant cone; independent coordinate scalings reduce it to $\mathbb R_+^4$ or to a product. Otherwise the boundary contains a point whose nonzero coordinates correspond to three noncollinear eigenvalue pairs, and the rank is three. The rotation case has rank drop only on the two coordinate projective lines. Thus, outside the already listed homogeneous and product cases, a two-dimensional boundary cannot be contained in the rank-drop set. At a rank-three boundary point, invariance makes the corresponding orbit relatively open in the boundary.
\end{proof}

\begin{lemma}[Independent scaling gives a product]\label{lem:product}
Suppose $\R^4=U\oplus\R e$ is invariant and the effective algebra contains an element represented as $diag(0,0,0,1)$ in a suitable basis. Then every proper invariant cone with nonempty interior is, up to replacing $e$ by $-e$, equivalent to a product
\[
 K_U\times\Rp.
\]
\end{lemma}

\begin{proof}
For $(u,se)\in K^\circ$, the relative scaling gives $(u,\e^tse)$. Letting $t\to-\infty$ gives $(u,0)\in K$. Multiplying by $\e^{-t}>0$ and letting $t\to+\infty$ gives $(0,se)\in K$. Let $K_U$ be the projection of $K$ to $U$ and let $K_e$ be its projection to $\R e$. The two limits show that $(u,0)\in K$ for every $u\in K_U$ and $(0,se)\in K$ for every $se\in K_e$. Since a cone is closed under addition, $(u,se)=(u,0)+(0,se)\in K$, so $K_U\times K_e\subset K$. The opposite inclusion follows from the definition of the projections. Pointedness makes $K_e$ a ray, which proves the statement.
\end{proof}

\begin{lemma}[Independent 2D scaling gives a product]\label{lem:product2D}
Suppose $\R^4 = V \oplus W$ is a direct sum of two $2$-dimensional invariant subspaces, $\dim V = 2$ and $\dim W = 2$, and the effective algebra contains an element represented as $\diag(1,1,0,0)$ in a suitable basis. Then every proper invariant cone $K \subset \R^4$ with nonempty interior is equivalent to a product of two $2$-dimensional proper cones
\[
 K = K_V \times K_W,
\]
where $K_V \subset V$ and $K_W \subset W$ are the projections of $K$ onto $V$ and $W$, respectively.
\end{lemma}

\begin{proof}
Let $(v, w) \in K^\circ$ with $v \in V$ and $w \in W$. The flow generated by the matrix $S = \diag(1,1,0,0)$ yields $e^{tS}(v, w) = (\e^t v, w) \in K^\circ$ for all $t \in \R$. Taking the limit as $t \to -\infty$, closedness of $K$ implies
\[
 \lim_{t \to -\infty} (\e^t v, w) = (0, w) \in K.
\]
Similarly, multiplying by $\e^{-t} > 0$ and taking the limit as $t \to +\infty$, we obtain
\[
 \lim_{t \to +\infty} \e^{-t} e^{tS}(v, w) = \lim_{t \to +\infty} (v, \e^{-t} w) = (v, 0) \in K.
\]
Let $K_V \subset V$ and $K_W \subset W$ be the projections of $K$ onto $V$ and $W$, respectively. The two limit arguments show that $(v, 0) \in K$ for every $v \in K_V$ and $(0, w) \in K$ for every $w \in K_W$. Since a convex cone is closed under vector addition,
\[
 (v, w) = (v, 0) + (0, w) \in K \implies K_V \times K_W \subset K.
\]
The reverse inclusion $K \subset K_V \times K_W$ holds trivially by the definition of projections. Therefore, $K = K_V \times K_W$. Since $K$ is proper, $K_V$ and $K_W$ are proper $2$-dimensional cones in $V$ and $W$, which completes the proof.
\end{proof}

%% file: sections/appendix_abelian.tex
\section{Complete analysis of the abelian case}\label{app:abelian}

Throughout this appendix $K\subset\R^4$ is a proper convex cone,
\[
 \Span(I,A,B)\subset\aut(K),\qquad [A,B]=0,
\]
and the action of $A,B$ on rays is effective and generically two-dimensional. In order to obtain a canonical form of the generators, we may add scalar matrices (multiples of $I$) to $A,B$, replace them by another basis of their span modulo $I$, and conjugate them simultaneously by a real invertible matrix, corresponding to a change of coordinates in $\mathbb R^4$.

\subsection{The invariant subspace decomposition} \label{sec:ProofTheorem3}

Here we provide a proof of Theorem \ref{thm:weights-main}.

\begin{proof}
Over $\mathbb C$, the commuting matrices can be triangularized simultaneously. Grouping equal pairs of diagonal entries gives common generalized eigenspaces. Complex conjugation either fixes such a space or joins it with its conjugate, producing the real decomposition~\eqref{eq:partitions-main}. 

A two-dimensional real space is either a direct sum of two complex-conjugate 1-dimensional eigenspaces, in which case we get the second form, or gives rise to a real Jordan block, in which case we get the first form. 

If the decomposition follows the pattern $3+1$, then the span of $I,A,B$ contains at least one non-zero nilpotent matrix. This matrix cannot have a Jordan block of size 2 by the nilpotent obstruction Lemma \ref{lem:nilpotent-obstruction}. Hence on the three-dimensional real subspace, the nilpotent matrix must be a Jordan block $J_3$. Solving $J_3X=XJ_3$ gives $X=aI+bJ_3+cJ_3^2$. 

This completes the proof of Theorem \ref{thm:weights-main}.
\end{proof}

\begin{proposition}[Exhaustion of the abelian case]\label{prop:abelian-full}
The abelian equation produces only the homogeneous and product cones in Tables~\ref{tab:classification-known}--\ref{tab:classification-lorentz}, the asymmetric power cones $K^{p,q,r}_{\alpha,\beta}$, the norm-power cones $K^{2,2}_\lambda$, the exponential-power cones $K_{\mathrm{pow},a}$, and the exponential-parabolic cone $K_{\mathrm{par}}$.
\end{proposition}

The proposition is proven in the subsequent sections, which consider the dimensional patterns of the invariant generalized eigenspace decomposition one by one.

\subsection{Four one-dimensional common eigenspaces} \label{sec:B2}

If all invariant subspaces are 1-dimensional, then both $A,B$ are diagonal. Let the pair of eigenvalues of $A,B$ corresponding to the basis vector $e_i$ be represented by
$\omega_i\in\R^2$. In order for $I,A,B$ to be linearly independent, we demand that the affine span of the four points $\omega_i$ equals $\R^2$. Then there is a unique affine relation
\begin{equation}\label{eq:weight-relation-app}
 \sum_{i=1}^4c_i=0,
 \qquad
 \sum_{i=1}^4c_i\omega_i=0,
\end{equation}
with coefficients $c_i$ determined up to a nonzero factor. It follows that the 0-homogeneous function $f(x) = \prod_{i=1}^4 |x_i|^{c_i}$ is an invariant of the group action. In each open orthant of $\mathbb R^4$, a level surface of this function determines a complete orbit.

Suppose that $c_i = 0$ for some $i$. Then $f$ does not depend on $x_i$, and the level surface is a direct product of a 2-dimensional surface with a ray $\mathbb R_{++}$. Since $K$ is closed and does not contain a line, in this case $x_i$ must either be always nonnegative or non-positive on $K$. Hence $K$ itself is a product cone, representable as a direct product of $\mathbb R_+$ with a 3-dimensional semi-homogeneous cone.

We therefore assume that all $c_i$ are nonzero. In particular, all four $\omega_i$ are distinct. We shall now consider subcases according to the signs of the $c_i$. The case when one coefficient has opposite sign to the other three has been already considered in Section \ref{sec:5}. Below we shall consider the remaining case when two coefficients are negative and two are positive.

%\subsubsection{Two coefficients of each sign} \label{sec:4abelian2plus2minus}

Let the coefficients $c_i$ of the affine relation \eqref{eq:weight-relation-app} be $-b,a,-d,c$ with $a,b,c,d > 0$, $a+c = b+d$. Then the generators can be chosen in the form
\[
 A=\diag(a,b,0,0),\qquad B=\diag(0,0,c,d),
\]
and the invariant can be chosen as
\begin{equation}\label{eq:two-two-diagonal-invariant-app}
 f=\frac{y^{a/(a+c)}w^{c/(a+c)}}{x^{b/(b+d)}z^{d/(b+d)}} =: \frac{g(y,w)}{h(x,z)}.
\end{equation}

Suppose that $K \neq \mathbb{R}_+^4$. Without loss of generality let $v \in \partial K \cap \operatorname{int}(\mathbb{R}_+^4)$ and $f(v) = \gamma$.
However, this level surface is neither convex nor concave. To exclude the first, choose two non-proportional positive points $(x_1,z_1)$ and $(x_2,z_2)$ with
$h(x_1,z_1)=h(x_2,z_2)=1$. Concavity of $h$ gives
\[
 h\left(\frac{x_1+x_2}{2},\frac{z_1+z_2}{2}\right)>1.
\]
Fix $(y,w)$ with $g(y,w)=\gamma$. The two resulting four-dimensional points $(x_k,y,z_k,w)$ lie on $f=\gamma$, while their midpoint has $f<\gamma$. The set $\{f\ge \gamma\}$ is therefore not convex. Interchanging numerator and denominator excludes $\{f\le \gamma\}$. We arrive at a contradiction with the convexity of $K$.

\subsection{One two-dimensional common subspace and two lines}

Write the coordinates as two line coordinates and one two-dimensional coordinate plane whenever this makes the final cone notation simpler.

\subsubsection{A rotation on the two-dimensional plane}\label{sec:abelian211rotation}

Each element of the $\Span(I,A,B)$ is represented in the form

\begin{equation}
 \begin{pmatrix}a&0&0&0\\0&b&0&0\\0&0&c&-d\\0&0&d&c\end{pmatrix}
\end{equation}

If $d \neq 0$, then it is impossible for $c > \max(a, b)$ or $c <  \min(a, b)$. Otherwise, by writing out the explicit orbit of some point and letting $t$ tend to $+\infty$ and $-\infty$ respectively, the cone would necessarily contain an entire plane. These two inequalities define opposite quadrants of 4-dimensional space, it means the three-dimensional subspace $\Span(I,A,B)$ must pass through the line $a = c = b$. 

Therefore we can normalize generator $B$ to the pure rotation, then normalize $A$, so that they would take the following form:
\begin{equation}\label{eq:norm-matrices-app}
 A=\begin{pmatrix}1&0&0&0\\0&0&0&0\\0&0&\lambda&0\\0&0&0&\lambda\end{pmatrix},
 \qquad
 B=\begin{pmatrix}0&0&0&0\\0&0&0&0\\0&0&0&-1\\0&0&1&0\end{pmatrix}
\end{equation}

Where $0 \le \lambda \le 1$ for the same reason as in the previous paragraph.

\begin{lemma}\label{lem:normPowerCone}
Let $K$ be a proper cone whose algebra corresponds to the given case. Then it is linearly isomorphic to the $K^{2,2}_\lambda$.
\end{lemma}
\begin{proof}
 All points in K must have the same sign for the $x$ coordinate, otherwise the cone will contain a line because of flow generated by $A$. For the same reason, using $I - A$, the sign of the $y$ coordinate is also the same everywhere in $K$.

 By changing coordinates, we make $x,y > 0$ in $K^\circ$. Let us prove that there is a point with $x>0, y>0, z>0$ in $\partial K$. If $(x_0, y_0, z_0, w_0) \in K^\circ$, then the ray $\{(x_0, y_0,z, 0) \mid z \geq z_0\}$ must intersect $\partial K$, otherwise $K$ will contain the affine plane $\{(x_0,y_0,z,w) \mid z,w \in \mathbb R\}$.

 The action of the flows generated by $I,A,B$ has the invariant
 \begin{equation}\label{eq:norm-invariant-app}
 f=\frac{\sqrt{z^2+w^2}}{x^\lambda y^{1-\lambda}}.
\end{equation}

Then the level set of $f$ with value $q_0 = \frac{z_0}{x_0^{\lambda}y_0^{1-\lambda}}$ at the boundary point $(x_0,y_0,z_0,0)$ with $x_0,y_0,z_0 > 0$ will be part of the boundary. The level sets with values in $[0,q_0)$ are in the interior, because the points $(x_0,y_0,z,w)$ with $z^2 + w^2 < z_0^2$ are in the convex hull of boundary points. On the other hand, the level sets with larger values are outside. 

By a scaling of $z,w$ the boundary value of the invariant can be made equal to 1.
\end{proof}

\subsubsection{A shear on the two-dimensional plane}\label{sec:B3.2}

Place the two-dimensional space in the $(x,z)$ coordinates and the two lines in the $y,w$ coordinates, so that $V = \Span\{e_1, e_3\}$. Suppose that there is no element $A$ of the algebra with distinct elements $A_{22} \not= A_{44}$. Then the nilpotent matrix $E_{13}$ is an element of the algebra, which contradicts Lemma \ref{lem:nilpotent-obstruction}. After scalar parts are removed, the general effective pair has the form
\begin{equation}\label{eq:shear-general-matrices-app}
 A_0=\begin{pmatrix}0&0&a&0\\0&1&0&0\\0&0&0&0\\0&0&0&0\end{pmatrix},
 \qquad
 B_0=\begin{pmatrix}0&0&b&0\\0&0&0&0\\0&0&0&0\\0&0&0&1\end{pmatrix}.
\end{equation}
 The signs of $y$ and $w$ are fixed in the interior of a proper invariant cone: the $A_0$-flow as $s\to+\infty$ isolates the $y$-ray after normalization by $\e^{-s}$, while the $B_0$-flow as $t\to+\infty$ isolates the $w$-ray. Coordinate reflections let us take $y,w>0$.
 
 Let us prove that while preserving the form of generators, we can move to a coordinate system where $(0,1,1,1) \in \partial K$.
 For some point $(x_0,y_0,z_0,w_0) \in \partial K$ we can make $z_0 = 1$ by changing coordinate $z$. Then, using $A, B$ flows we can make $y_0 = 1, w_0 = 1$. If $x_0 \neq 0$ yet, then by changing coordinates $y' = ye^{x_0/a}$ and using A-flow we set $x_0$ to zero.
 
 Write the orbit through $(0,1,1,1)$ as
\begin{equation}\label{eq:shear-orbit-app}
 \{(as+bt,\e^s,1,\e^t):s,t\in\R\}.
\end{equation}
In the section $z=1$, it is the graph
\begin{equation}\label{eq:shear-graph-app}
 x=a\log y+b\log w.
\end{equation}
The Hessian is
\[
 D^2(a\log y+b\log w)
 =\begin{pmatrix}-a/y^2&0\\0&-b/w^2\end{pmatrix}.
\]
By Lemma~\ref{lem:graph-convexity}, a convex side exists only when $a,b$ have the same sign. A zero coefficient gives a product by lemma \ref{lem:product}. Replacing $x$ by $-x$ makes $a,b>0$. If the interior contained points with both signs of $z$, apply the two flows with $s=t=-T$ and let $T\to+\infty$. The $y$- and $w$-coordinates tend to zero, while $x$ equals $-(a+b)Tz+O(1)$. Division by $T$ would put both $e_1$ and $-e_1$ in $K$. Hence the sign of $z$ is fixed as well, and after a final reflection we take $z>0$. Conjugating in the $(x,z)$ plane normalizes the coefficient of the shear in $A_0$ to one, and replacing $B_0$ by $B_0-(b/a)A_0$ gives
\begin{equation}\label{eq:Kpow-matrices-app}
 A_{\mathrm{pow}}=\begin{pmatrix}0&0&1&0\\0&1&0&0\\0&0&0&0\\0&0&0&0\end{pmatrix},
 \qquad
 B_{\mathrm{pow}}=\begin{pmatrix}0&0&0&0\\0&-\mu&0&0\\0&0&0&0\\0&0&0&1\end{pmatrix},
 \qquad \mu=\frac ba>0.
\end{equation}

\begin{lemma}\label{lem:K_pow_mu}
Let K be a proper semihomogeneous cone whose algebra corresponds to the given case. Then it is linearly isomorphic to one of the cones of the family $\{ K_{\mathrm{pow}, \mu} \mid 0 < \mu \leq 1 \}$
\end{lemma}
\begin{proof}
The invariant of $\Span\{I, A_{\mathrm{pow}}, B_{\mathrm{pow}}\}$ is 
\begin{equation}\label{eq:Kpow-invariant-app}
 f=\log(y/z)-x/z+\mu\log(w/z).
\end{equation}

Choose a point $(0,1,1,1) \in \partial K$. The orbit will be a level set of $f$, given by
\[
 y=z\e^{x/z}(w/z)^{-\mu}.
\]

On $z=1$, the boundary function $\phi(x,w)=\e^xw^{-\mu}$ satisfies
\begin{equation}\label{eq:Kpow-Hessian-app}
 D^2\phi=\phi\begin{pmatrix}1&-\mu/w\\-\mu/w&\mu(\mu+1)/w^2\end{pmatrix},
 \qquad
 \det D^2\phi=\frac{\mu \phi^2}{w^2}>0.
\end{equation}
Thus only the epigraph is convex, which corresponds to $K_{\mathrm{pow},\mu}$. The linear map
\begin{equation}\label{eq:a-inversion-app}
 (x,y,z,w)\longmapsto(x/\mu,w,z,y)
\end{equation}
sends $K_{\mathrm{pow},\mu}$ to $K_{\mathrm{pow},1/\mu}$, so $0<\mu\le1$ is a complete range.

\end{proof}

\subsection{Two two-dimensional common subspaces}\label{subsec:mixed22}
Combine the rotation $\begin{psmallmatrix}0&-1\\1&0\end{psmallmatrix}$ on $(x,y)$ with the shear $\begin{psmallmatrix}0&1\\0&0\end{psmallmatrix}$ on $(z,w)$. If the two coefficients are independent, the matrices are
\begin{equation}\label{eq:mixed-independent-matrices-app}
 A=\begin{pmatrix}a&-1&0&0\\1&a&0&0\\0&0&0&0\\0&0&0&0\end{pmatrix},
 \qquad
 B=\begin{pmatrix}b&0&0&0\\0&b&0&0\\0&0&0&1\\0&0&0&0\end{pmatrix}.
\end{equation}
Choose $(x_0, y_0, z_0, w_0) \in \partial K$, where $w_0 \neq 0$. Using $I$ flow we can make $w_0 = 1$. Using $A$ and $B$ flows we can make $x_0 = 1, y_0 = 0, z_0 = 0$. The orbit through $(1,0,0,1)$ is
\begin{equation}\label{eq:mixed-independent-orbit-app}
 \left\{(\e^{a\theta+bs}\cos\theta,
 \e^{a\theta+bs}\sin\theta,s,1):\theta,s\in\R\right\}.
\end{equation}
Near $(\theta,s)=(0,0)$, use $(y,z)$ as local coordinates and write the first coordinate as $x=f(y,z)$. Direct differentiation gives
\begin{equation}\label{eq:mixed-independent-Hessian-app}
 D^2f(0,0)=\begin{pmatrix}-(1+a^2)&0\\0&b^2\end{pmatrix}.
\end{equation}
The two signs exclude both sides by Lemma~\ref{lem:graph-convexity}.

If the two coefficients are proportional, a complementary relative scaling gives
\begin{equation}\label{eq:mixed-proportional-matrices-app}
 A=\begin{pmatrix}0&-1&0&0\\1&0&0&0\\0&0&0&1\\0&0&0&0\end{pmatrix},
 \qquad
 B=\diag(1,1,0,0).
\end{equation}
Choose $(x_0, y_0, z_0, w_0) \in \partial K$, where $w_0 \neq 0$, then using $I$ make $w_0 = 1$. Using $A,B$ flows make $x_0 = 1, y_0 = 0, z_0 = 0$.The orbit through $(1,0,0,1)$ is
\begin{equation}\label{eq:mixed-proportional-orbit-app}
 \{(\e^s\cos\theta,\e^s\sin\theta,\theta,1):s,\theta\in\R\}.
\end{equation}
At every point with $\sin\theta\ne0$, this is the graph
\[
 x=y\cot z.
\]
Its Hessian satisfies
\begin{equation}\label{eq:mixed-proportional-Hessian-app}
 \det D^2(y\cot z)=-\csc^4z<0.
\end{equation}
Again neither side is convex.

Finally, let both two-dimensional restrictions contain a nonzero rotation. If there exists an element of the algebra with distinct diagonal elements, then the dynamics on one two-dimensional subspace dominates the other under the corresponding flow, and $K$ has to contain a two-dimensional linear subspace. If all elements of the algebra have a constant diagonal, then the trace-less subalgebra is $\mathfrak{so}(2) \times \mathfrak{so}(2)$. The orbits of the generated subgroup are tori which contain the origin in the interior of their convex hull.

Therefore the pattern $2+2$ gives no cone.

\subsection{A three-dimensional common subspace and one line} \label{sec:B5}

Choose a nonzero matrix in the kernel of the relative line character. By Lemma~\ref{lem:nilpotent-obstruction}, its restriction to the three-dimensional space is the matrix in~\eqref{eq:three-matrix-app}. Solving the commutation equation and changing generators gives
\begin{equation}\label{eq:three-one-matrices-app}
 A=\begin{pmatrix}0&1&0&0\\0&0&1&0\\0&0&0&0\\0&0&0&0\end{pmatrix},
 \qquad
 B=\begin{pmatrix}0&0&c&0\\0&0&0&0\\0&0&0&0\\0&0&0&1\end{pmatrix}.
\end{equation}

By a coordinate change $(x,z) \mapsto (\alpha x,\alpha^{-1}z)$ we may scale $c$ by a positive multiplier, so we need to consider only the cases $c = 0,\pm\frac12$.

For $c=0$, Lemma~\ref{lem:product} splits off the fourth coordinate. On the remaining three coordinates, the flow of $\begin{psmallmatrix}0&1&0\\0&0&1\\0&0&0\end{psmallmatrix}$ preserves $2xz-y^2$. In the section $z=1$ its orbits are the parabolas $x=y^2/2+c_0$. The only convex invariant side with nonempty interior is the epigraph $x\ge y^2/2$, whose conic closure is $L_3$. Hence the product is $L_3\times\Rp$. 

Let now $c \not= 0$. We wish to find a boundary orbit containing the point $(0,0,1,1)^T$. We may assume $w \geq 0$, otherwise the action of $B$ creates $\pm e_4 \in K$. Likewise $z \geq 0$, otherwise the action of $A$ creates $\pm e_1 \in K$. If $z = 0$, then also $y = 0$, otherwise the action of $A$ creates $\pm e_1 \in K$. Hence there exists a boundary point with $z,w > 0$. The action of $A$ can be used to force $y = 0$, while the action of $B$ can be used to force $x = 0$. Finally, the scaling of $w$ achieves $z = w$ on the boundary.

For $c\ne0$, the orbit through $(0,0,1,1)$ is
\begin{equation}\label{eq:three-one-orbit-app}
 \{(s^2/2+c\log w,s,1,w):s\in\R,\ w>0\}.
\end{equation}
In the section $z=1$, the graph is
\begin{equation}\label{eq:three-one-graph-app}
 x=\frac{y^2}{2}+c\log w,
 \qquad
 D^2f=\begin{pmatrix}1&0\\0&-c/w^2\end{pmatrix}.
\end{equation}
For $c>0$, the two signs exclude both sides. For $c<0$, the epigraph is convex. Rescale to $c=-1/2$ and apply the invertible map
\begin{equation}\label{eq:Kpar-linear-map-app}
 (x,y,z,w)\longmapsto(-2x,w,z,y).
\end{equation}
The inequality becomes
\[
 y\ge z\exp\left(\frac xz+\left(\frac wz\right)^2\right),
\]
which is $K_{\mathrm{par}}$. The normalized-flow limits show that the signs of $y$ and $z$ cannot change inside a proper invariant cone; after coordinate reflections they are positive. Thus this epigraph is the complete invariant cone, not merely one local boundary patch. In the new coordinates the matrices are
\begin{equation}\label{eq:Kpar-matrices-app}
 A_{\mathrm{par}}=\begin{pmatrix}0&0&1&0\\0&1&0&0\\0&0&0&0\\0&0&0&0\end{pmatrix},
 \qquad
 B_{\mathrm{par}}=\begin{pmatrix}0&0&0&-2\\0&0&0&0\\0&0&0&0\\0&0&1&0\end{pmatrix}.
\end{equation}
The boundary orbit through $(0,1,1,0)$ is
\begin{equation}\label{eq:Kpar-orbit-app}
 \{(s-t^2,\e^s,1,t):s,t\in\R\}.
\end{equation}
The Hessian of $\e^{x+w^2}$ is
\begin{equation}\label{eq:Kpar-Hessian-app}
 D^2\e^{x+w^2}=\e^{x+w^2}
 \begin{pmatrix}1&2w\\2w&4w^2+2\end{pmatrix},
 \qquad \det=2\e^{2x+2w^2}>0.
\end{equation}

\subsection{One four-dimensional common subspace}\label{sec:B6}
Assume first that the common eigenvalue pair is real. After scalar subtraction, the effective algebra is a two-dimensional space of commuting nilpotent matrices. Every nonzero matrix must have a three-dimensional block and a one-dimensional zero block. Choose
\begin{equation}\label{eq:four-A-app}
 A=\begin{pmatrix}0&1&0&0\\0&0&0&1\\0&0&0&0\\0&0&0&0\end{pmatrix}.
\end{equation}
Solving $CA=AC$ for a general matrix $C=(c_{ij})$ gives
\begin{equation}\label{eq:four-commuting-matrix-app}
 C=\begin{pmatrix}
 a&b&c&d\\0&a&0&b\\0&0&e&f\\0&0&0&a
 \end{pmatrix}.
\end{equation}
After scalar subtraction and adding a multiple of $A$, the second nilpotent generator is
\[
 B=cE_{13}+dE_{14}+fE_{34}.
\]
Conjugation by a matrix $D = I_4 + \frac{d}{c}E_{34}$ commuting with $A$ removes $d$. Direct multiplication then gives
\begin{equation}\label{eq:four-square-app}
 (sA+tB)^2=(s^2+cf\,t^2)E_{14}.
\end{equation}
If $cf=0$, a nonzero square-zero matrix appears. If $cf<0$, choose $(s,t)\ne(0,0)$ with $s^2+cf\,t^2=0$. Both cases contradict Lemma~\ref{lem:nilpotent-obstruction}. Hence $cf>0$, and scaling coordinates and generators gives
\begin{equation}\label{eq:four-final-matrices-app}
 A=\begin{pmatrix}0&1&0&0\\0&0&0&1\\0&0&0&0\\0&0&0&0\end{pmatrix},
 \qquad
 B=\begin{pmatrix}0&0&1&0\\0&0&0&0\\0&0&0&1\\0&0&0&0\end{pmatrix}.
\end{equation}
These matrices preserve $w$ and $2xw-y^2-z^2$. If the interior contained points with both signs of $w$, the quadratic terms in the two flows, divided by the square of the flow parameter, would put both $e_1$ and $-e_1$ in $K$. Thus the sign of $w$ is fixed; replace $w$ by $-w$ if necessary. In the section $w=1$, the orbit with invariant value $q$ is
\begin{equation}\label{eq:four-orbit-app}
 \{((s^2+t^2+q)/2,s,t,1):s,t\in\R\}.
\end{equation}
The only convex side is the epigraph
\[
 2xw\ge y^2+z^2,\qquad w\ge0.
\]
The change $x\mapsto x-qw/2$ normalizes the endpoint to zero. This is $L_4$.

If the common eigenvalue pair is nonreal, a nonzero generalized part has, in real coordinates, the square-zero matrix
\[
 \begin{pmatrix}0&0&1&0\\0&0&0&1\\0&0&0&0\\0&0&0&0\end{pmatrix},
\]
which is excluded. Without this matrix, both generators are combinations of the identity and one simultaneous rotation, so the action on rays has dimension at most one. Together with Lemma~\ref{lem:regular-boundary-direct}, which guarantees that a nonproduct cone has one of the regular boundary orbits analyzed above, this proves Proposition~\ref{prop:abelian-full}.

%% file: sections/appendix_affine.tex
\section{Complete analysis of the equation $[A,B]=B$}\label{app:affine}

In this section we consider the case when the principal orbits of $\Aut(K)$ have dimension 3 and the 3-dimensional connected subgroup $G \subset \Aut(K)$ with orbits of dimension 3 is not abelian. Throughout this section we assume
\[
 \Span(I,A,B) = \mathfrak g \subset \aut(K),\qquad [A,B]=B,
\]
where $K\subset\R^4$ is proper and the action of $A,B$ on interior rays of $K$ is effective.

\subsection{The matrix $B$ and the full commutator system} \label{sec:commutatorSystem}

\begin{lemma}\label{lem:B-nilpotent}
The matrix $B$ is nilpotent.
\end{lemma}

\begin{proof}
Induction gives $[A,B^k]=kB^k$. Taking traces yields
$\tr(B^k)=0$ for $k=1,\ldots,4$. Newton's identities then show that every coefficient of the characteristic polynomial except the leading one is zero. Hence every eigenvalue of $B$ is zero.
\end{proof}

The nonzero (the dimension of the principal orbit excludes $B=0$) nilpotent Jordan types in dimension four have block dimensions
\[
 4,\qquad3+1,\qquad2+2,\qquad2+1+1.
\]
Lemma~\ref{lem:nilpotent-obstruction} excludes the first, third, and fourth. Therefore, possibly after a change of coordinates, we may assume
\begin{equation}\label{eq:B-affine-app}
 B=\begin{pmatrix}0&1&0&0\\0&0&1&0\\0&0&0&0\\0&0&0&0\end{pmatrix}.
\end{equation}

\begin{lemma}[Entrywise solution of the commutator equation]\label{lem:commutator-solution}
Let $A=(a_{ij})$. The equation $AB-BA=B$ is equivalent to Table~\ref{tab:commutator}. Its complete solution is
\begin{equation}\label{eq:commutator-general-app}
 A=
 \begin{pmatrix}
 \mu+2&b&c&d\\
 0&\mu+1&b&0\\
 0&0&\mu&0\\
 0&0&e&\mu+\gamma
 \end{pmatrix}.
\end{equation}
\end{lemma}

\begin{table}[h]
\caption{All nonzero scalar equations in $AB-BA-B=0$.}
\label{tab:commutator}
\centering
\small
\setlength{\tabcolsep}{5pt}
\renewcommand{\arraystretch}{1.12}
\begin{tabular}{@{}llll@{}}
\toprule
$(1,1)$: $a_{21}=0$ & $(1,2)$: $a_{11}-a_{22}=1$ & $(1,3)$: $a_{12}=a_{23}$ & $(1,4)$: $a_{24}=0$\\
$(2,1)$: $a_{31}=0$ & $(2,2)$: $a_{21}=a_{32}$ & $(2,3)$: $a_{22}-a_{33}=1$ & $(2,4)$: $a_{34}=0$\\
$(3,2)$: $a_{31}=0$ & $(3,3)$: $a_{32}=0$ & $(4,2)$: $a_{41}=0$ & $(4,3)$: $a_{42}=0$\\
\bottomrule
\end{tabular}
\end{table}

\begin{proof}
The two nonzero columns of $B$ and the two nonzero rows of $B$ give exactly the equations in Table~\ref{tab:commutator}. They imply
\[
 a_{21}=a_{31}=a_{32}=a_{24}=a_{34}=a_{41}=a_{42}=0,
 \quad a_{12}=a_{23},
 \quad a_{11}=a_{33}+2,
 \quad a_{22}=a_{33}+1.
\]
The free entries are $a_{33},a_{12},a_{44},a_{13},a_{14},a_{43}$, which gives~\eqref{eq:commutator-general-app}. The coefficient matrix has rank ten and nullity six.
\end{proof}

The parameters $\mu,b$ in \eqref{eq:commutator-general-app} can be removed by adding multiples of $I$ and $B$ to $A$. To simplify the remaining entries, use
\begin{equation}\label{eq:commuting-conjugation-app}
 C(u,v,q)=
 \begin{pmatrix}0&0&u&v\\0&0&0&0\\0&0&0&0\\0&0&q&0\end{pmatrix},
 \qquad C^3=0,
\end{equation}
which commutes with $B$. Direct calculation of
$\e^{C}A\e^{-C}$ gives
\begin{align}
 d'&=d+(\gamma-2)v,\label{eq:dprime}\\
 e'&=e-\gamma q,\label{eq:eprime}\\
 c'&=c-dq+ev+(1-\gamma)vq-2u.\label{eq:cprime}
\end{align}
If $\gamma\ne0,2$, choose $v,q,u$ in this order to make $d'=e'=c'=0$. At $\gamma=0$, only $e$ may remain; at $\gamma=2$, only $d$ may remain. This proves Proposition~\ref{prop:affine-main}.

Below we consider the emerged subcases.

\subsection{Case $c = d = e = 0$}\label{sec:C2}

Set
\begin{equation}\label{eq:zero-tail-matrices-app}
 A_\gamma=\begin{pmatrix}2&0&0&0\\0&1&0&0\\0&0&0&0\\0&0&0&\gamma\end{pmatrix},
 \qquad
 B_N=\begin{pmatrix}0&1&0&0\\0&0&1&0\\0&0&0&0\\0&0&0&0\end{pmatrix}.
\end{equation}
The generated flows are
\begin{align}
 \e^{tB_N}(x,y,z,w)&=(x+ty+t^2z/2,y+tz,z,w),\label{eq:BL-flow-app}\\
 \e^{sA_\gamma}(x,y,z,w)&=(\e^{2s}x,\e^sy,z,\e^{\gamma s}w).\label{eq:Agamma-flow-app}
\end{align}
The first preserves $\Delta=2xz-y^2$, and the second sends $\Delta$ to $\e^{2s}\Delta$.

If $z > 0$ at some interior point, then $e^{tB_N}$ sends this point in the direction of $e_1$ for $t \to \pm\infty$. Likewise, the orbit under the action of $e^{tB_N}$ of a point with $z < 0$ tends to $-e_1$. Hence $z$ cannot have opposite signs in $K$, and we may assume $z > 0$ in $K^{\circ}$ without loss of generality.

\begin{lemma}\label{lem:affine-chambers}
Let $K$ be a full-dimensional proper cone invariant under $I,A_\gamma,B_N$, on which $z \geq 0$. After replacing $K$ by its image under a coordinate reflection, the coordinate $z$ is positive on $K^\circ$. The projection of $K^\circ$ to the $(x,y,z)$-space has, in the section $z=1$, exactly one of the two forms
\[
 \{(x,y) \mid 2x-y^2>0\},
 \qquad\text{or}\qquad
 \R^2.
\]
\end{lemma}

\begin{proof}
Let $P$ be the projection of $K^\circ$ to the first three coordinates and put
$S=P\cap\{z=1\}$. Linear projection is open, so $S$ is open and convex. The ensemble of the two base flows in~\eqref{eq:BL-flow-app}--\eqref{eq:Agamma-flow-app} is transitive on each of the sets
\[
 2x-y^2>0,
 \qquad 2x-y^2=0,
 \qquad 2x-y^2<0.
\]
The last open set is not convex: $(0,1)$ and $(0,-1)$ belong to it, whereas their midpoint does not. Therefore an invariant open convex set is either the first set or all of $\R^2$.
\end{proof}

\subsubsection{Case $S = \mathbb R^2$}

\begin{lemma}
    Assume the notations of the previous lemma and its proof and suppose that the set $S$ equals $\mathbb R^2$. Then for $\gamma>2$ the only nonproduct possibility is $K^N_{\mathrm{ext},\gamma}$, while $\gamma=2$ gives the product $L_3\times\Rp$ and the cone $K^N_{2,1,1}$. All other values of $\gamma$ are impossible.
\end{lemma}

\begin{proof}
    Assume now that $S=\R^2$. Convexity and invariance imply that, after possibly replacing $w$ by $-w$, the section $z = 1$ of $K$ is the epigraph of a convex function $w = \phi(x,y)$ defined on $\mathbb R^2$ and whose value depends only on $r=2x-y^2$. The variable $w$ cannot be bounded from above and from below simultaneously, because then by homogeneity of the bounds we would get $w = 0$ at $r = 0$. The scaling by $A_\gamma$ gives
\[
 \phi\bigl(\e^{2s}r\bigr)=\e^{\gamma s}\phi(r).
\]
Thus, with $p=\gamma/2$,
\begin{equation}\label{eq:full-base-endpoint}
 \phi=c_+r^p\quad(r>0),
 \qquad
 \phi=c_-(-r)^p\quad(r<0)
\end{equation}
for some constants $c_{\pm}$, with the evident interpretation at $\gamma=0$. By convexity $\phi$ is continuous at $r=0$.

If $\gamma < 2$, then the action of $A_{\gamma}$ sends a generic point in the limit $s \to +\infty$ in the direction $\pm e_1$, depending on the sign of $x$ at the initial point. Hence $x$ cannot assume both signs on $K$, and $S = \mathbb R^2$ is ruled out. 

Let $\gamma>2$. At $(x,y) = (1,0)$, the function $\phi$ is strictly convex in the direction of $x$ and strictly concave in the direction of $y$ if $c_+ > 0$. Hence $c_+=0$. Then $c_- > 0$, otherwise we get again a line in $K$. Since the function $\phi = \max(0,(y^2-2x)^p)$ is convex, we indeed obtain an invariant convex cone, which after a normalization of the coordinate $w$ turns out to be the cone $K_{\mathrm{ext},\gamma}^N$ from Table \ref{tab:classification-lorentz}.

For $\gamma=2$ we get $p = 1$, and convexity of $\phi$ implies that $c_+ \leq 0$, $c_- \geq 0$, and $c_+ + c_- \geq 0$. If $c_+ + c_- = 0$, then $\phi = c_- \cdot (y^2 - 2x)$, and $K$ contains the line $\{ (x,0,1,-2c_-x) \mid x \in \mathbb R \}$, contradicting pointed-ness. If $c_+ = 0$, $c_- > 0$, then $\phi = \max(0,c_- \cdot (y^2-2x))$, and the cone $K$ is determined by the inequalities $z,w \geq 0$, $(w+2c_-x)z \geq c_-y^2$. It is not hard to see that it is isomorphic to $L_3 \times \mathbb R_+$.

Let us consider the case $c_- > -c_+ > 0$. The cone $K$ is determined by the inequalities
\[ z \geq 0,\quad wz \geq |c_{\pm}| \cdot (y^2 - 2xz).
\]
Replacing the coordinates $w,x$ by $w' = \frac{w + 2c_-x}{c_-}$, $x' = \frac{w - 2c_+x}{-c_+}$, this transforms to $w'z \geq y^2$, $x'z \geq y^2$, $z \geq 0$. It is not hard to see that this cone is isomorphic to the cone $K^N_{2,1,1}$.
\end{proof}

Put
\begin{equation}\label{eq:exterior-function-app}
 \widetilde H_\delta(x,y,z)
 =z^{1-\delta}(y^2-2xz)_+^{\delta/2}.
\end{equation}
The function $y^2-2x$ is convex on the section $z=1$, its positive part is convex, and $r\mapsto r^{\delta/2}$ is convex and increasing. Therefore
$(y^2-2x)_+^{\delta/2}$ is convex. Formula~\eqref{eq:exterior-function-app} is its perspective and is convex and homogeneous. Consequently
\begin{equation}\label{eq:KL-ext-app}
 K^N_{\mathrm{ext},\delta}
 =\cl\{(x,y,z,w)\in\R^4:z>0,\ w>\widetilde H_\delta(x,y,z)\}
\end{equation}
is a convex cone. On the open exterior chamber its complete invariant is
\begin{equation}\label{eq:KL-ext-invariant-app}
 q^{\mathrm{ext}}_\delta
 =\frac{w}{z^{1-\delta}(-\Delta)^{\delta/2}}.
\end{equation}
A lower endpoint gives the epigraph. An upper endpoint is a nonconvex hypograph, and scaling removes the positive lower constant. Thus~\eqref{eq:KL-ext-app} is the only nonproduct cone in the exterior chamber.

\subsubsection{Case $S = \{(x,y) \mid 2x > y^2\}$}

Now suppose that $S = \{(x,y) \mid 2x > y^2\}$, in particular, $\Delta > 0$ on $K^{\circ}$. For given $(x,y,z) = (x,y,1)$ with $2x > y^2$, the variable $w$ takes values in an interval (which may be finite or infinite) which is not equal to the real line, otherwise $K$ is not pointed. Hence there exist boundary points in $K$ with $z > 0$, $\Delta > 0$. The orbit invariant is given by
\begin{equation}\label{eq:L-invariant-app}
 q_\gamma=\frac{w}{H_\gamma},
 \qquad H_\gamma=\Delta^{\gamma/2}z^{1-\gamma}.
\end{equation}

In every orbit on which $z>0$, $\Delta>0$, there exists a point with $(x,y,z)=(1/2,0,1)$. Indeed, $B_N$ sends every point with fixed $(z,\Delta,w)$ to another such point with $y=0$, and $A_\gamma,I$ can be used for scaling of $x,z$. A level orbit is then given by
\begin{equation}\label{eq:L-orbit-app}
 \left\{\rho \cdot \left(
 \frac{\e^{2s}(1+t^2)}2,\e^st,1,c\e^{\gamma s}
 \right):\rho>0,\ s,t\in\R\right\},
\end{equation}
where $c$ is the value of the invariant $q_{\gamma}$ on the orbit. In particular, $\partial K$ contains orbits of this form.

Let us verify the convexity of the boundary of $K$. Let $(1/2,0,1,w_0)$ be a point on the boundary orbit. Then the level surface of $H_{\gamma}$ through $(1/2,0,1)$ is the boundary of $K$ intersected with the hyperplane $w = w_0$.

At $(x,y,z)=(1/2,0,1)$, direct differentiation of $H_{\gamma}$ gives the Hessian
\begin{equation}\label{eq:Hessian-main}
 D^2H_\gamma=
 \begin{pmatrix}
 \gamma(\gamma-2)&0&\gamma(2-\gamma)/2\\
 0&-\gamma&0\\
 \gamma(2-\gamma)/2&0&\gamma(\gamma-2)/4
 \end{pmatrix}.
\end{equation}

Its kernel is the radial vector $(1/2,0,1)$. The signs of the two other eigenvalues equal those of $-\gamma$ and $\gamma(\gamma-2)$ and cannot be opposite. This yields $\gamma^2(\gamma-2) \leq 0$, or equivalently $\gamma \leq 2$. Because the group generated by~\eqref{eq:BL-flow-app}--\eqref{eq:Agamma-flow-app} is transitive on the rays of $\cL^\circ$, the same signs hold everywhere. Below we shall consider subcases according to the value of $\gamma$.

\textbf{$\gamma = 2$, $\Delta > 0$: product cone and $K^N_{2,1,1}$}

Then $H_{\gamma} = \frac{2xz-y^2}{z}$. This function is concave, and $w$ is bounded from above by a nonnegative multiple and from below by a non-positive multiple of $H_{\gamma}$. These multiples may also be infinite. This leads to the cones $K^N_{2,\alpha,\beta}$, $(\alpha,\beta) \in {\cal E}$. The transformation $x \mapsto x + \mu w$ preserves the general form of the inequalities $-\alpha\Delta \leq wz \leq \beta\Delta$ but leads to $\alpha \mapsto \frac{\alpha}{1+2\mu\alpha}$, $\beta \mapsto \frac{1}{1-2\mu\beta}$. Hence all cones with both $\alpha,\beta$ positive are equivalent to $K^N_{2,1,1}$. If $\alpha = 1$, $\beta = +\infty$, then the coordinate change $w \gets \frac{w+2x}{2}$ leads again to $K^N_{2,1,1}$. If $\alpha = 0$, $\beta = 1$, then $x \gets 2x - w$ leads to the product cone $L_3 \times \mathbb R_+$. If $\alpha = 0$, $\beta = +\infty$, then we directly get $L_3 \times \mathbb R_+$.

\textbf{$0<\gamma<2$, $\Delta > 0$: Lorentz-power cones}

The function $H_\gamma$ is concave and homogeneous of degree one. Therefore
\[
 w\le\beta H_\gamma
\]
defines a convex set exactly for $\beta\ge0$, and
\[
 w\ge-\alpha H_\gamma
\]
defines a convex set exactly for $\alpha\ge0$. For every fixed $(x,y,z) = (x,y,1)$ with $2x > y^2$ the value of $w$ must vary in an interval, not equal to the real line and containing 0. We obtain
\begin{equation}\label{eq:KL-app}
 K^N_{\gamma,\alpha,\beta}
 =\cl\{(x,y,z,w) \mid x,z>0,\Delta>0,-\alpha H_\gamma<w<\beta H_\gamma\},
 \qquad(\alpha,\beta)\in\EndPairs.
\end{equation}
The base $\cL$ is pointed, and the continuous closure over its vertex contains at most one transverse ray, so the full cone is pointed. At $\gamma=1$ and $(\alpha,\beta)=(1,1)$ the inequality is
\[
 2xz\ge y^2+w^2,
\]
which defines $L_4$.

\textbf{$\gamma=0, \Delta > 0$: product cone and $K^N_{0,1,1}$}

In this case $H_{\gamma} = z$. On the hyperplane $z = 1$ the variable $w$ may be bounded by constants from below and from above. If $w$ is unbounded, then $K$ contains a line. If $w$ is bounded from one side, then $K$ is a product cone. If $w$ lies in a compact interval, then this interval may be shifted to $[-1,1]$, and $K$ is defined by the inequalities $z \geq 0$, $2xz \geq y^2$, $-z \leq w \leq z$, which corresponds to $K^N_{0,1,1}$. 

\textbf{$\gamma<0, \Delta > 0$: the Lorentz epigraph}

Now $H_\gamma$ is convex and positive. Thus the epigraph
\begin{equation}\label{eq:KL-epi-app}
 K^N_{\mathrm{mix},\gamma}
 =\cl\{(x,y,z,w) \mid x,z>0,\Delta>0,w>H_\gamma\}
\end{equation}
is convex. Any other positive multiple is removed by scaling $w$. A finite upper endpoint would require a hypograph of a strictly convex function and is not convex; two finite endpoints contain such an upper inequality. Reflection in $w$ treats the opposite orientation. Hence~\eqref{eq:KL-epi-app} is the only nonproduct cone in this case. 

\subsection{The two logarithmic matrices}\label{sec:C3}

\subsubsection{The cone $K^N_{\log}$}

At $\gamma=0$, normalize the nonzero entry $e$ in~\eqref{eq:commutator-general-app} by a scaling of $w$ to obtain
\begin{equation}\label{eq:Alog-app}
 A_{\log}=\begin{pmatrix}2&0&0&0\\0&1&0&0\\0&0&0&0\\0&0&2&0\end{pmatrix},
 \qquad B_N=\begin{pmatrix}0&1&0&0\\0&0&1&0\\0&0&0&0\\0&0&0&0\end{pmatrix}.
\end{equation}
The flow generated by $A_{\log}$ is given by
\[
 (x,y,z,w)\longmapsto(\e^{2s}x,\e^sy,z,w+2sz).
\]
This yields $\Delta\mapsto\e^{2s}\Delta$, and $\Delta$ has constant sign on every orbit. As in the case $c = d = e = 0$, the coordinate $z$ cannot take both signs in $K$, otherwise the action of $B_N$ leads to $\pm e_1 \in K$, contradicting pointed-ness. Hence we assume $z > 0$ on $K^{\circ}$ without loss of generality. As in Lemma \ref{lem:affine-chambers}, we may distinguish the two cases $S = \mathbb R^2$ and $S = \{(x,y) \mid 2x-y^2 > 0\}$. However, the action of $A_{\log}$ sends a generic point in the direction $\pm e_1$ in dependence on the sign of $x$ at the initial point. Hence $x$ also cannot assume both signs in $K$, and the case $S = \mathbb R^2$ has to be ruled out.

We then may assume $x,z,\Delta > 0$ on $K^{\circ}$. The orbit invariant is given by
\begin{equation}\label{eq:log-invariant-app}
 q_{\log} = \frac{w}{z}-\log\frac{\Delta}{z^2}.
\end{equation}
The level set $q_{\log}=q_0$ determines
\[
 w=z \cdot (q_0 + \log(\Delta/z^2))
\]
as a function of $x,y,z$ on $\{ (x,y,z) \mid z,\Delta > 0 \}$. This function is concave for all $q_0 \in \mathbb R$, so $w$ can only be upper bounded by a level set on $K$.

By a change of coordinates $w \gets w - q_0z$ we may assume that this level set is given by $q_0 = 0$, which yields the cone $K^N_{\log}$.

%The convex side $q_{\log}<0$ is
%\begin{equation}\label{eq:Klog-epigraph-app}
% x>\frac{y^2}{2z}+\frac z2\e^{w/z},
%\end{equation}
%which is the epigraph of a sum of two convex perspectives. 

\subsubsection{The cone $K^N_{\log,2}$}

At $\gamma=2$, normalize the nonzero entry $d$ by a scaling of $w$ to
\begin{equation}\label{eq:Alog2-app}
 A_{\log,2}=\begin{pmatrix}2&0&0&2\\0&1&0&0\\0&0&0&0\\0&0&0&2\end{pmatrix},
 \qquad B_N=\begin{pmatrix}0&1&0&0\\0&0&1&0\\0&0&0&0\\0&0&0&0\end{pmatrix}.
\end{equation}
As in the previous case, we may assume $z > 0$ on $K^{\circ}$ without loss of generality by considering the flow generated by $B_N$.

The flow generated by $A_{\log,2}$ is given by
\[
 (x,y,z,w)\longmapsto
 (\e^{2s}(x+2sw),\e^sy,z,\e^{2s}w).
\]
For $s \to +\infty$ this flow drives a generic point in the direction $\pm e_1$ in dependence of the initial sign of $w$. Hence $w$ cannot assume both signs on $K$ and we may assume $w > 0$ on $K^{\circ}$ ($z$ can flip the sign independently of $w$).

Hence
\[
 \Delta\longmapsto\e^{2s}(\Delta+4swz),
\]
and the orbit invariant is
\begin{equation}\label{eq:Klog2-invariant-app}
 q_{\log,2}=\frac{\Delta}{2wz}-\log\frac wz.
\end{equation}

The level set $q_{\log,2} = q_0$ determines
\[ x = w \cdot \left(q_0 + \frac{y^2}{2wz} + \log\frac{w}{z}\right)
\]
as a function of $y,z,w$. This function is convex for all $q_0$, and can hence only lower bound $x$ on $K$ to assure convexity.

By a coordinate change $x \gets x - q_0w$ we may assume that the boundary of $K$ is given by the level set $q_{\log,2} = 0$. This yields
\begin{equation}\label{eq:Klog2-epigraph-app}
 x>\frac{y^2}{2z}+w\log(w/z),
 \qquad z,w>0,
\end{equation}
corresponding to the cone $K^N_{\log,2}$.

%% file: sections/appendix_verification.tex
\section{Normalizer actions and uniqueness of parameters}

\begin{proposition}[Parameter identifications]\label{prop:parameters}
The parameter identifications stated in the third column of Tables~\ref{tab:classification-known} and~\ref{tab:classification-lorentz} are the only linear isomorphisms between parameter values. Distinctly named nonproduct families are not linearly isomorphic, except for the homogeneous overlaps explicitly removed in the tables.
\end{proposition}

\begin{proof}
A linear isomorphism between two proper cones $K_1 \cong K_2$ conjugates their full automorphism Lie algebras: $\aut(K_2) = M \aut(K_1) M^{-1}$ for some $M \in \mathrm{GL}_4(\mathbb{R})$. The dimensions of the common invariant subspaces and the displayed non-scalar matrices intrinsically separate the abelian families: $K^{p,q,r}_{\alpha,\beta}$ has four common $1$-dimensional eigenlines; $K^{2,2}_\lambda$ has one $2$-dimensional plane on which the displayed rotation matrix acts and two eigenlines; $K_{\mathrm{pow},a}$ has the displayed shear matrix on a $2$-dimensional plane and two eigenlines; $K_{\mathrm{par}}$ has a $3$-dimensional common invariant space and one eigenline. Product cones contain an independent relative scaling as in Lemma~\ref{lem:product}. The nonabelian families are defined by the relation $[A,B]=B$.

For $K^{p,q,r}_{\alpha,\beta}$, the four common eigenlines are intrinsic. A change of the two projective generators applies an affine transformation to the four eigenvalue pairs, and a linear normalizer only permutes the four eigenlines. The unique affine relation has coefficients
\[
 (1/p, 1/q, 1/r, -1),
\]
up to a common factor and permutation. Hence, the unordered triple $\{p,q,r\}$ is invariant.

For $K^{2,2}_\lambda$, the $2$-dimensional plane and the two eigenlines are intrinsic. The only nontrivial normalizer action interchanges the two eigenlines, sending $\lambda$ to $1-\lambda$. For $K_{\mathrm{pow},a}$, the coefficient of the shear and the two relative line weights determine the positive ratio $a$, up to interchange of the two line coordinates, which realizes $a \leftrightarrow a^{-1}$. The cone $K_{\mathrm{par}}$ has no parameter.

For the algebra generated by $A_\gamma$ and $B_N$, every algebra automorphism sends $B_N$ to a nonzero multiple of itself and $A_\gamma$ to $A_\gamma + c B_N$ modulo scalars. The flag of submodules
\[
 0 \subset \ker B_N \cap \operatorname{im} B_N^2 \subset \operatorname{im} B_N^2 \subset \operatorname{im} B_N \subset \mathbb{R}^4
\]
and the $1$-dimensional complementary eigenspace are intrinsic to the action of $B_N$. The successive weights of $A_\gamma$ on the $3$-dimensional block are $2, 1, 0$, while its weight on the complementary line is $\gamma$. Hence, $\gamma$ is invariant. For $\gamma \neq 0, 1, 2$, solving the matrices that commute with $B_N$ shows that $q_\gamma = w / H_\gamma$ changes only by a nonzero scalar and possibly a sign, so the normalized endpoint pair $(\alpha, \beta) \in \mathcal{E}$ is complete.

At $\gamma = 1$, the full algebra is $\mathbb{R} I_4 \oplus \mathfrak{sl}_2(\mathbb{R})$ and the representation has the intrinsic decomposition
\[
 \mathbb{R}^4 = \operatorname{Sym}^2(\mathbb{R}^2) \oplus \mathbb{R}.
\]
Every normalizer preserves these two spaces. The discriminant $\Delta = 2xz - y^2$ is, up to a positive scalar on the Lorentz chamber, the unique quadratic relative invariant on the first space, and $w$ is multiplied by a nonzero scalar on the second. Therefore, $w / \sqrt{\Delta}$ changes only by a nonzero scalar and possibly a sign. The interval remains a complete invariant.

At $\gamma = 0$ and $\gamma = 2$, the affine and fractional-linear formulas in Section~\ref{app:affine} give exactly the endpoint identifications. The nonzero exceptional entries in $A_{\log}$ ($a_{43} = 2$) and $A_{\log,2}$ ($a_{14} = 2$) cannot be removed by any matrix commuting with $B_N$, so the logarithmic cones $K^N_{\log}$ and $K^N_{\log,2}$ do not belong to the zero-tail family (where $c=d=e=0$). Furthermore, they are not isomorphic to one another: the boundary face of $K^N_{\log}$ at $z = 0$ is $2$-dimensional, while the boundary face of $K^N_{\log,2}$ at $z = 0$ is $1$-dimensional. The ranges $\gamma < 0$ and $\gamma > 2$ are separated by the invariant value of the complementary weight $\gamma$, so the two epigraph families $K^N_{\mathrm{mix},\gamma}$ and $K^N_{\mathrm{ext},\gamma}$ do not overlap.

The product parameters are complete by the $3$-dimensional classification in \citet{HildebrandAffine2014}. This proves the proposition.
\end{proof}
\newpage
\section{Affine compact sections of some cones} \label{sec:naming}
\begin{figure}[htbp]
    \centering
    \begin{minipage}[c]{0.60\textwidth}
        \centering
        \includegraphics[width=\textwidth]{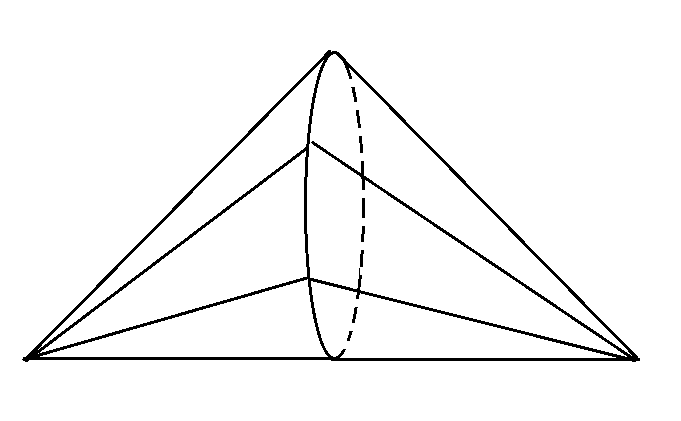}
    \end{minipage}
    \hfill
    \begin{minipage}[c]{0.36\textwidth}
        \centering
        \includegraphics[width=\textwidth]{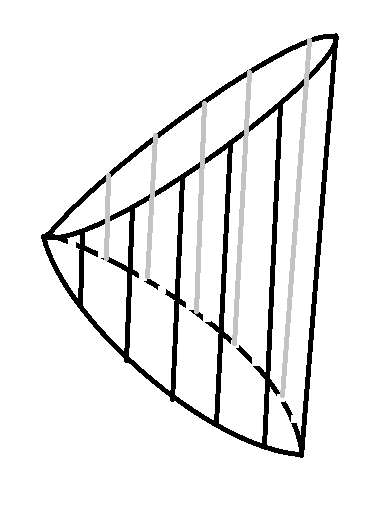}
    \end{minipage}
    
    \caption{The affine compact section of the cone $K^N_{2,1,1}$ resembles a croissant (left), that of $K^N_{0,1,1}$ a cylindrical wedge (right).}
    \label{fig:geometry1}
\end{figure}

In this section we provide pictures of the affine compact sections of some of the singleton cones, which do not belong to a continuous family. The sections of the cones $K^N_{2,1,1}$ and $K^N_{0,1,1}$ are depicted on Fig.~\ref{fig:geometry1}, explaining their naming. These two cones are dual to each other.

\begin{figure}[htbp]
    \centering
    \begin{minipage}[c]{0.55\textwidth}
        \centering
        \includegraphics[width=\textwidth]{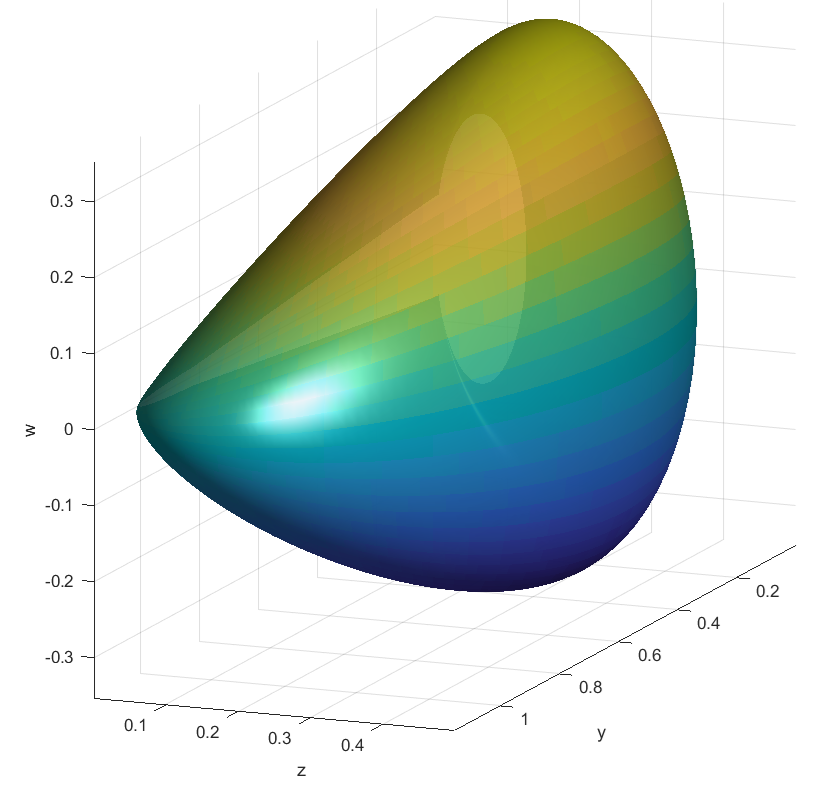}
    \end{minipage}
    \hfill
    \begin{minipage}[c]{0.40\textwidth}
        \centering
        \includegraphics[width=\textwidth]{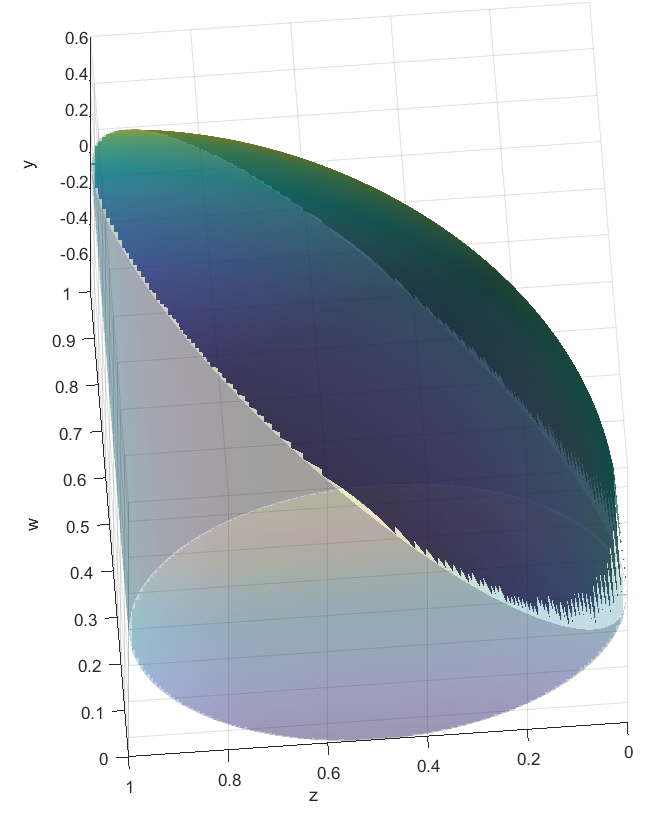}
    \end{minipage}
    
    \caption{The affine compact section $y+z=1+x$ of the exponential-quadratic cone (left) and the section $y+z = 1$ of the exponential-parabolic cone (right). The intersection of the plane $y = 0$ with the section of the exponential-quadratic cone is a point given by $z = w = 0$, that of the plane $z = 0$ a segment of the line $w = y = 0$. Numerically these intersections resemble an ellipse and a triangle due to exponential approach to the supporting planes.}
    \label{fig:geometry2}
\end{figure}

On Fig.~\ref{fig:geometry2} we depict the sections of the exponential-quadratic and the exponential-parabolic cones.

